\documentclass{siamart220329}

\usepackage{amssymb}
\usepackage{graphicx}
\usepackage{epstopdf} 
\usepackage{hyperref}    
\usepackage{subcaption}
\usepackage{amssymb,latexsym,amsmath}
\usepackage{xcolor}
\usepackage{enumerate}
\usepackage{booktabs}
\usepackage{MnSymbol,bbding,pifont}

\usepackage{tikz,tikz-cd}
\usetikzlibrary{hobby}
\usetikzlibrary{arrows}
\usetikzlibrary{positioning}
\usetikzlibrary{shapes,snakes}
\usetikzlibrary{cd}

\def\R{\mathbb{R}}
\def\C{\mathbb{C}}
\def\N{\mathbb{N}}

\DeclareMathOperator{\rank}{rank}

\DeclareMathOperator{\adj}{adj}

\DeclareMathOperator{\tr}{trace}

\DeclareMathOperator{\spann}{span}

\newtheorem{remark}[theorem]{{\sc Remark}}
\newtheorem{example}[theorem]{Example}
\newtheorem{problem}[theorem]{Problem}

\newcommand{\hide}[1]{}

\numberwithin{equation}{section}

\begin{document}

\title{A Riemannian optimization method to compute the nearest singular pencil\thanks{Submitted to the editors on March 5th, 2024}}

\author{
Froil\'{a}n Dopico\thanks{Departamento de Matem\'aticas, Universidad Carlos III de Madrid, Avda. Universidad 30, 28911 Legan\'{e}s, Spain (\email{dopico@math.uc3m.es}). Supported by { the ``Agencia Estatal de Investigaci\'on'' of Spain through grants  PID2019-106362GB-I00 MCIN/AEI/10.13039/501100011033 and RED2022-134176-T,} and by the Madrid Government (Comunidad de Madrid-Spain) under the Multiannual Agreement with UC3M in the line of Excellence of University Professors (EPUC3M23), and in the context of the V PRICIT (Regional Programme of Research and Technological Innovation).}
\and Vanni Noferini\thanks{Aalto University, Department of Mathematics and Systems Analysis, P.O. Box 11100, FI-00076, Aalto, Finland (\email{vanni.noferini@aalto.fi}). Supported by an Academy of Finland grant (Suomen Akatemian p\"{a}\"{a}t\"{o}s 331230).}
\and Lauri Nyman\thanks{Corresponding author. Aalto University, Department of Mathematics and Systems Analysis, P.O. Box 11100, FI-00076, Aalto, Finland (\email{lauri.s.nyman@aalto.fi}).}}

\date{}

\maketitle

\begin{abstract}
Given a square pencil $A+ \lambda B$, where $A$ and $B$ are $n\times n$ complex {(resp. real)} matrices, we consider the problem of finding the singular {  complex (resp. real)} pencil nearest to it in the Frobenius distance. This problem is known to be very difficult, and the few algorithms available in the literature can only deal efficiently with pencils of very small size. We show that the problem is equivalent to minimizing a certain objective function { $f$} over the Riemannian manifold $SU(n) \times SU(n)$ {(resp. $SO(n) \times SO(n)$ if the nearest real singular pencil is sought)}, where $SU(n)$ denotes the special unitary group {(resp. $SO(n)$ denotes the special orthogonal group)}. This novel perspective is based on the generalized Schur form of pencils, and yields { competitive numerical methods}, by pairing it with { algorithms} capable of doing optimization on { Riemannian manifolds. We propose one algorithm that directly minimizes the (almost everywhere, but not everywhere, differentiable) function $f$, as well as a smoothed alternative and a third algorithm that is smooth and can also solve the problem} of finding a nearest singular pencil with a specified minimal index. We provide numerical experiments that show that the resulting { methods allow} us to deal with pencils of much larger size than alternative techniques, yielding candidate minimizers of comparable or better quality. In the course of our analysis, we also obtain a number of new theoretical results related to the generalized Schur form of a (regular or singular) square pencil and to the minimal index of a singular square pencil whose nullity is $1$.
\end{abstract}

\begin{keywords}
regular matrix pencil, singular matrix pencil, optimization on matrix manifolds, nearness matrix problems, generalized Schur form, minimal indices
\end{keywords}

\begin{AMS}
65F99, 65K10, 15A18, 15A22
\end{AMS}

\section{Introduction}

Let $A,B \in \C^{n \times n}$ be a pair of square matrices over the complex field. Linear polynomial matrices of the form $A + \lambda B$ are called \emph{matrix pencils} \cite{gant59}. Matrix pencils are important in a variety of applications; namely, it is often required to solve the generalized eigenvalue problem, which consists of finding the values $\mu \in \C \cup \{\infty \}$ such that $\rank (A+\mu B) < \rank (A+\lambda B)$, where we conventionally agree $A+\infty B :=B$ and $\rank (A+\lambda B)$ is the normal rank of the pencil, that is, the largest size of a minor of $A+ \lambda B$ which is not the identically zero polynomial in $\lambda$. Sometimes, it is also of interest to compute other quantities such as eigenvectors or, if appropriate, minimal indices and minimal bases. Generalized eigenvalue problems arise in many applications including for example discretization of PDEs, control, differential algebraic equations, and linearization of more general nonlinear eigenvalue problems. See, for instance, \cite{kunkel-merhmann06,vand81} and the references therein.  The square matrix pencil $A + \lambda B$ is called \emph{regular} if $\det (A + \lambda B) \not \equiv 0$ and \emph{singular} otherwise.  It is immediate from this definition that singular $n \times n$ pencils are a proper Zariski closed subset of the set of all $n \times n$ pencils, and in particular they are therefore nowhere dense.

On the other hand, there is a remarkable difference between regular and singular pencils when dealing with the associated eigenvalue problems \cite{vand81}.  Recently some progress has been made to advance both the theory \cite{DN20,LN20} and the numerical solution \cite{HMP,KS23,LN20} of singular generalized eigenvalue problems. Nevertheless, in the regular case, both the spectral theory and the design and analysis of numerical methods simplify considerably. Furthermore, for certain applications regularity is crucial. For instance, in the context of systems of linear algebraic-differential equations, the existence of a unique solution to the original problem is only guaranteed when the associated pencil is regular \cite{gg,KV15,kunkel-merhmann06}. This motivates the problem, given a regular pencil, to compute the distance to the nearest singular pencil of the same size as well as the singular pencil that minimizes this distance, which was posed in \cite{bhm}. In the literature, some attention has also been received by some variants of the problem, such as for example certain structured versions of the nearest singular pencil problem \cite{ps22} or the extension to polynomial matrices of possibly higher degrees \cite{bora,gg}.

The quest for the nearest singular pencil was first launched by Byers, He and Mehrmann \cite{bhm}, who collected several lower and upper bounds and proposed several characterizations of the distance to singularity but could not provide an exact solution except in selected very special cases. Later, several numerical methods to approximate a solution were proposed \cite{bora,ghl,glm,KV15,MU14}. However, the tasks of computing the distance to singularity and/or a nearest singular pencil were soon recognised to be extremely difficult: the number of local minima seems to increase fast with the size of the pencil (see Subsection \ref{sec:multiple_starting} for some experiments),
making it hard to find the global optimum. 

Moreover, the existing numerical algorithms are not particularly efficient; for example each iteration of the algorithm in \cite{ghl} has an asymptotic complexity of $O(n^{12})$ flops for an $n \times n$ input pencil. While there is no explicit complexity analysis for the algorithms in \cite{bora,glm,MU14}, the heuristics is not promising when $n$ is large. It was observed in \cite{KV15} that the method of \cite{MU14} can exhibit very slow convergence and is only practical for very small sizes. The algorithm in \cite{glm} has been perhaps the best available method for some years, but no efficient implementation is available and in \cite{glm} the algorithm could only be tested against extremely small examples ($n \leq 8$). In the very recent paper \cite{bora}, significant progress was achieved as a method based on ideas already mentioned in \cite{bhm} was implemented efficiently, and not only for pencils but generally for polynomial matrices of any degree: To quote the authors, the resulting algorithm is efficient ``for problems that are not large" \cite[Section 10]{bora}, and again \cite{bora} only reports tests on extremely small inputs ($n \leq 8$).
 
Generally, the existing methods rely on either (i) ODE-based techniques \cite{gg,glm}, or (ii) structured perturbations of (potentially very large) block Toeplitz matrices containing the coefficients of the pencil as blocks \cite{bhm,ghl,MU14}, or (iii) other related structured optimization problems also involving matrices containing the coefficients of the pencil as blocks \cite{bora}. In this paper, we take instead a completely different approach by following the lead of Noferini and Poloni, who recently applied Riemannian optimization techniques to compute the matrix nearest to a given one and having all the eigenvalues in a prescribed region $\Omega$ \cite{np}. Evidence was given in \cite{np} that, for the nearest $\Omega$-stable matrix problem, the Riemannian optimization approach dramatically outperformed other numerical methods for input matrices of small to moderate size (up to a few hundreds). In \cite{np} the nearest stable matrix problem was recast as an optimization problem on the manifold of unitary (or orthogonal) matrices. In this paper, we crucially show that several of the main ideas of \cite{np} can be extended to the nearest singular pencil problem. Namely, we prove that also this problem is equivalent to { minimizing a certain objective function} on a Riemannian manifold: in this case, the manifold is the Cartesian product of the manifold of unitary matrices times itself. This strategy allows us to devise an algorithm that (a) utterly outshines its predecessors in terms of speed, with the only exception of pencils of very small size where the method of \cite{bora} can be faster, (b) often provides an improvement in the quality of the computed minima, and (c) allows us, for the first time, to efficiently tackle inputs whose size is not extremely small (the largest size for which we ran experiments is $n = 200$, see Subsection \ref{sec:running_time}).
{
However, the objective function $f$ to be minimized in our main algorithm is non-smooth because one of its addends is the minimum of $n$ smooth functions. For this reason, while no practical issues were observed during extensive experiments, a rigorous proof of convergence for every input is currently unavailable. Motivated by this, we also propose two smooth alternatives. The first is a smoothed version of the direct algorithm, which replaces the minimum addend in $f$ with the Boltzmann softmin operator \cite{asadi17}. Numerical tests suggest that this approach is slower and yields minimizers of similar quality as the direct approach. The second alternative solves $n$ smooth minimization subproblems, each obtained by replacing the minimum addend in $f$ with the $i$th argument of the minimum, and selects the final solution(s) corresponding to the index (or indices) $i$ that gives the smallest $f_i$. This latter approach is slower but provides additional information, as each subproblem solves the restricted problem of finding a nearest singular pencil with a specified minimal index \cite[Ch. XII]{gant59}. MATLAB codes for our direct nonsmooth method and the smoothed versions are available at \href{https://github.com/NymanLauri/nearest-singular}{https://github.com/NymanLauri/nearest-singular}. In practice, we find that the nonsmooth method is the best option.

}

It is worth mentioning that, in some applications, the problem of finding the nearest \emph{real} singular pencil to a given real pencil is also of interest. { This real variant of the problem can be solved by minimizing exactly the same cost function as in the complex case but over the Cartesian product of the manifold of orthogonal matrices times itself. Moreover, the complex minimization algorithm over the Cartesian product of the manifold of unitary matrices times itself solves automatically the real variant of the problem for a given real pencil whenever a real starting point is provided.}
% \begin{vn}
%     Above, I removed reference to Sections. The description of the paper's structure is below, so those comments sounded redundant here.
% \end{vn}
% \begin{fd}
%     This is OK for me.
% \end{fd}

\hide{
%However, while the basic ideas are similar, certain key parts of the algorithm become much more technical in the real case due to the presence of $2\times 2$ diagonal blocks in the real generalized Schur form. To keep the length of the present article reasonable, we thus postpone the analysis of the real version of the method to future research, and we only focus on the complex case here.
}

The structure of the manuscript is as follows. In Section \ref{sec:problem}, we state the nearest singular pencil problem precisely and provide some initial observations. In Section \ref{sec:complex}, we show that the problem is equivalent to the minimization of a certain objective function on the manifold of ordered pairs of unitary matrices; consequently, we devise an algorithm capable of solving the latter optimization task. {Section \ref{sec:smoothalgorithm} describes a  smoothed variant of our main algorithm.}  In Section \ref{sec:third}, { we discuss how to solve the nearest singular pencil problem by solving $n$ smooth minimization subproblems. Each subproblem corresponds to the problem, discussed in Section \ref{sec:fixminindex},} 
of finding the nearest singular pencil having a prescribed value of its (right) minimal { index. It turns out} that a simple modification of our algorithm can solve this problem too by means of Riemannian optimization, and along the way we also prove some theoretical results on square singular pencils that might be interesting per se. { Sections \ref{sec:problem}-\ref{sec:fixminindex} deal with complex pencils, but in Section \ref{sec:real} we prove a result on the generalized Schur form of real singular pencils that allows us to extend all the results obtained for complex pencils to the problem of computing a singular real pencil nearest to a given real pencil.} In Section \ref{sec:imp} we briefly describe { the practical implementation of our algorithms}, and in Section \ref{sec:numerical} we report a number of numerical experiments to highlight the features of the proposed { methods}. We finally draw some conclusions in Section \ref{sec:conclusion}. Some technical results on Schur forms of regular pencils are given in Appendix \ref{sec:swapping}, { while Appendix \ref{sec:gradient_smooth} contains the explicit expressions of the Euclidean gradient and Hessian of the smoothed objective function of Section \ref{sec:smoothalgorithm}}.

% \begin{fd}
% Check what is decided on Appendix \ref{sec:gradient_smooth}.
% \end{fd}

% \begin{vn}
%     I would be OK to move Appendix \ref{sec:gradient_smooth} to the answer to referees, but I think Froilan (and maybe also Lauri?) disagree and prefer to keep it in the paper. My opinion is not strong on this, and I can be happy with your decision.
% \end{vn}

% \begin{fd}
% I prefer to keep the Appendix \ref{sec:gradient_smooth}. If I were a reader, I would like to see the appendix just to evaluate the difficulty of the approach and for having some possibility to check the correctness of the smoothed algorithm.  
% \end{fd}

\section{The nearest singular pencil problem}\label{sec:problem}

A pencil over $\C$ is a polynomial of degree at most $1$ whose coefficients are matrices (of the same size) with complex entries. We define $\C[\lambda]_1^{n\times n}$ to be the vector space of $n \times n$ pencils over $\C$, that is, the space  $\C[\lambda]_1^{n \times n}:=\{A+\lambda B$ $|$ $A,B \in \C^{n \times n}\}$. Clearly, $\C[\lambda]_1^{n \times n} \cong \C^{2n^2}$, and therefore $\C[\lambda]_1^{n \times n}$ can be equipped with the Euclidean norm $\| A + \lambda B \|_F$ whose square is defined as follows
\begin{equation}\label{eq:pencilnorm}
    \| A + \lambda B \|_F^2 := \sum_{i,j=1}^n (|A_{ij}|^2+|B_{ij}|^2) = \left\| \begin{bmatrix}
A\\
B
\end{bmatrix} \right\|_F^2 = \left\| \begin{bmatrix}
A&B
\end{bmatrix}  \right\|_F^2,
\end{equation} 
where $\|\cdot \|_F$ is the Frobenius matrix norm and, for $a \in \C$, $|a|$ is the modulus of $a$. The norm \eqref{eq:pencilnorm} induces a distance 
\begin{equation}\label{eq:pencildist}
    d(A+\lambda B,S+\lambda T):=\|(A-S)+\lambda(B-T) \|_F.
\end{equation}
Throughout the paper, we always refer to \eqref{eq:pencildist} when speaking of the distance between two pencils or of the singular pencil nearest to a given one. The problem of our interest can be formulated as follows:

\begin{problem}\label{problem}
Given $A+\lambda B \in \C[\lambda]_1^{n \times n}$, compute both the minimum and a minimizer of the distance function $d(A+\lambda B,S+\lambda T)$ \eqref{eq:pencildist} amongst all pencils $S+\lambda T \in \C[\lambda]_1^{n \times n}$ that are \emph{singular}, that is, $\det(S+\lambda T) \equiv 0$.
\end{problem}
For practical purposes, in the following it will be convenient to express Problem \ref{problem} as the minimization of the \emph{square} distance; the equivalence is clear since $x \mapsto x^2$ is increasing on $[0,\infty[$. More explicitly, we equivalently express the problem as the computation of
\begin{equation}\label{eq:riproblem}
    \min_{S+\lambda T \in \mathcal{S}_n} [d(A+\lambda B,S+\lambda T)]^2,
\end{equation}
as well as its argument minimum, where $\mathcal{S}_n$ denotes the subset of the singular pencils in $\C[\lambda]_1^{n \times n}${, that is, 
\[
\mathcal{S}_n := \{S+\lambda T \in \C[\lambda]_1^{n \times n} \ | \ \det(S+\lambda T) \equiv 0\}.
\]
}

\section{A Riemannian algorithm for the nearest complex singular pencil}\label{sec:complex}
\subsection{Reformulating the problem} \label{subsec:reformulating}

We denote now by $U(n)$ the set of $n \times n$ unitary matrices and by $SU(n)$ the special unitary group, that is, the set of $n \times n$ unitary matrices with determinant $1$. Our most basic tool is the generalized Schur form of a pencil.
\begin{lemma}[Stewart \cite{Stewart}] \label{lemma:schur}
For any pair $A,B \in \C^{n \times n}$ there exist $Q,Z \in U(n)$ such that $QAZ$ and $QBZ$ are both upper triangular.
\end{lemma}
The upper triangular pencil $QAZ + \lambda QBZ$ given by the matrices in Lemma \ref{lemma:schur} is called a generalized Schur form of $A + \lambda B$. It is easy to see that each pencil has infinitely many generalized Schur forms. {Indeed, a Schur form can be multiplied by arbitrary unitary diagonal matrices on both sides to yield another Schur form; moreover, diagonal elements can be swapped  (see Corollary \ref{cora2} for the regular case).}

Next, we aim to extend the ideas first developed in \cite{np} to Problem \ref{problem}. To this goal, we first observe that it is easy to characterize singular triangular pencils.
\begin{lemma} \label{lemma:singular}
An upper triangular square pencil $A + \lambda B$ is singular if and only if it has at least one zero diagonal element.
\end{lemma} 
\begin{proof}
    It suffices to observe that the determinant of a triangular pencil is the product of the diagonal elements of the pencil.
\end{proof}
    
Lemma~\ref{lemma:singular} allows us to find the nearest singular \textit{upper triangular} pencil in a straightforward manner. This is outlined in Proposition~\ref{prop:nearest_tri}.
\begin{proposition} \label{prop:nearest_tri}
Let $A+\lambda B \in \C[\lambda]_1^{n \times n}$. Let $k$ be any index such that $|A_{kk}|^2+|B_{kk}|^2$ is minimal over the diagonal entries of $A+ \lambda B$.  An upper triangular singular pencil nearest to $A+\lambda B$ is $\mathcal{P}(A) + \lambda \mathcal{P}(B)$ where
\[ \mathcal{P}(A)_{ij}=\begin{cases}A_{ij} \ &\mathrm{if} \ i<j \ \mathrm{or} \  i=j\neq k;\\
0 \ &\mathrm{otherwise}; \end{cases} \qquad \mathcal{P}(B)_{ij}=\begin{cases}B_{ij} \ &\mathrm{if} \ i<j \ \mathrm{or} \  i=j\neq k;\\
0 \ &\mathrm{otherwise}. \end{cases}\]
In particular, the squared distance to $A+\lambda B$ from $\mathcal{P}(A)+\lambda \mathcal{P}(B)$ is 
\begin{equation}\label{eq:objectivefunction} \mathcal{F}(A+\lambda B)= \sum_{i>j} (|A_{ij}|^2+|B_{ij}|^2) + \min_{1 \leq i \leq n} (|A_{ii}|^2+|B_{ii}|^2 ).\end{equation}
\end{proposition}
\begin{proof}
For upper triangularity, it is necessary that the strictly lower triangular part of the minimizer is zero. To minimize the distance, one cannot do any better than preserving the strictly upper triangular part of $A,B$ in the minimizer, which is certainly possible. Finally, it is necessary for the singularity to annihilate in the minimizer at least one diagonal element of $A+\lambda B$, and to annihilate just one is neccessary for minimizing the distance; the definition of $k$ makes sure that the distance is indeed minimized.
\end{proof}

Fix now $A+\lambda B \in \C[\lambda]_1^{n \times n}$. Proposition~\ref{prop:nearest_tri} can be used, with the help of Lemmas~\ref{lemma:schur} and \ref{lemma:singular}, to characterize the singular pencil $S+\lambda T \in \C[\lambda]_1^{n \times n}$ nearest to $A+\lambda B$.  To this goal, recall that $\mathcal{S}_n \subset \C[\lambda]_1^{n \times n}$ is the subset of singular pencils and denote by $\mathcal{T}_n \subset \C[\lambda]_1^{n \times n}$ the subset of singular upper triangular pencils. For any pair of unitary matrices $(Q,Z) \in U(n) \times U(n)$ define the function
\begin{equation}\label{eq:newobjfun}
     f(Q,Z) := \mathcal{F}(QAZ+\lambda QBZ) = [d(QAZ+\lambda QBZ,\mathcal{P}(QAZ)+\lambda \mathcal{P}(QBZ))]^2  
\end{equation}
where $\mathcal{F}$ and $\mathcal{P}$ are as in Proposition~\ref{prop:nearest_tri} and $d$ is as in \eqref{eq:pencildist}. It turns out that Problem \ref{problem} is equivalent to computing
\begin{equation}\label{eq:minnewobjfun}
    \min_{(Q,Z)\in U(n) \times U(n)} f(Q,Z) 
\end{equation} 
as well as an argument minimum $(Q_0,Z_0)$. Theorem \ref{thm:obj_function} below states this result more formally.

\begin{theorem} \label{thm:obj_function}
Let $A+\lambda B \in \C[\lambda]_1^{n \times n}$. Let $f(Q,Z)$ be the function on $U(n) \times U(n)$ defined by \eqref{eq:newobjfun}. Then:
\begin{enumerate}
    \item The optimization problems \eqref{eq:riproblem} and \eqref{eq:minnewobjfun} have the same minimum value;
    \item The pair of unitary matrices $(Q_0,Z_0)$ is a global (resp. local) minimizer for \eqref{eq:minnewobjfun} if and only if the pencil $Q_0^* \mathcal{P}(Q_0 A Z_0) Z_0^* + \lambda Q_0^* \mathcal{P}(Q_0 B Z_0) Z_0^*$ is a global (resp. local) minimizer for \eqref{eq:riproblem}.
\end{enumerate}
\end{theorem}
\begin{proof}
Taking a generalized Schur form of the optimization variables of \eqref{eq:riproblem}, say, $X=QSZ,Y=QTZ$ we obtain the equivalence
\[ \min_{S+\lambda T \in \mathcal{S}_n} \|(A-S)+\lambda(B-T)\|_F^2  = \min_{Q,Z \in U(n)} \min_{X+\lambda Y  \in \mathcal{T}_n}  \|(A-Q^*XZ^*)+\lambda(B-Q^*YZ^*) \|_F^2 \]
\[ = \min_{Q,Z \in U(n)} \min_{X+\lambda Y  \in \mathcal{T}_n}  \|(QAZ-X)+\lambda(QBZ-Y) \|_F^2 = \min_{Q,Z \in U(n)} \mathcal{F}(QAZ+\lambda QBZ),\]
having used Proposition \ref{prop:nearest_tri} in the last step. This equation proves item 1 immediately; and, using again Proposition \ref{prop:nearest_tri}, item 2 also follows.
\end{proof}

This yields a method: we can minimize the objective function $f(Q,Z)$ over $U(n) \times U(n)$ to compute the nearest singular pencil. Any algorithm capable of doing optimization on Riemannian manifolds can thus be employed, as long as we input the right manifold and objective function.

{ \begin{remark} \label{rem:realthmmain} The proof of Theorem \ref{thm:obj_function} only requires the upper triangular generalized Schur form of the {\em singular} complex pencils $S + \lambda T$. We anticipate that we will prove in Lemma \ref{lem:froilannew} that such upper triangular generalized Schur form also exists for {\em singular real} pencils via real orthogonal equivalence. Therefore, Theorem \ref{thm:obj_function}, as well as the other results in this section, remain valid for real pencils, mutatis mutandis. %Readers who are particularly interested on real pencils may find useful to read Section \ref{sec:real} now, before proceeding with the next developments. 
\end{remark}}

\begin{remark}
  Observe that, in the statement of Theorem \ref{thm:obj_function}, we can replace $U(n)$ by $SU(n)$. To see why, it suffices to consider the following refinement of Lemma \ref{lemma:schur}.
    \begin{lemma}\label{lem:SUnotU}
    In Lemma \ref{lemma:schur}, without loss of generality one can take $Q,Z \in SU(n)$.
\end{lemma}
\begin{proof}
      Let $Q',Z' \in U(n)$ be such that $Q'AZ',Q'BZ'$ are upper triangular, and define $Q=\left(I_{n-1} \oplus \frac{1}{\det Q'}\right) Q'$, $Z=Z'\left(I_{n-1} \oplus \frac{1}{\det Z'}\right)$. Clearly $Q,Z \in SU(n)$ and $QAZ,QBZ$ are also upper triangular.
\end{proof}
It follows that we may implement our method using $SU(n) \times SU(n)$ as our search space. In the implementation used for the numerical experiments in this paper, however, we have used the manifold $U(n) \times U(n)$, for the practical reason that $U(n)$, but not $SU(n)$, is built-in in the MATLAB toolbox Manopt \cite{manopt}. This might have led to suboptimal running times in the numerical tests presented in Section \ref{sec:numerical}.
\end{remark}

It is worth noting that, even when the problem is cast in the form \eqref{eq:minnewobjfun}, it remains extremely difficult to solve numerically. Indeed, the optimization problem \eqref{eq:minnewobjfun} is highly non-convex, and as such our method can only yield a \textit{locally} minimal solution to the problem. {Moreover, the objective function \eqref{eq:newobjfun} is not smooth due to the addend $\min_{1 \leq i \leq n} (|(QAZ)_{ii}|^2 + |(QBZ)_{ii}|^2)$, which might not be differentiable at points where the minimum is attained at more than one value of $i$. Since \eqref{eq:newobjfun} is Lipschitz, this only happens on a measure-zero set of values of $(Q,Z)$; moreover, as we show below, a point of non-differentiability cannot coincide with a local minimum.

    \begin{theorem}\label{thm:lasttheorem}
        For every input $A + \lambda B$, the objective function $f(Q,Z)$ in \eqref{eq:newobjfun} is real-differentiable at all its local minima.
    \end{theorem}
    \begin{proof}
Note that \eqref{eq:newobjfun} is equivalent to \eqref{eq:fqz}. The only possible source of non-smoothness in the objective function arises from the min function, which can be non-differentiable at points $x_0=(Q_0,Z_0)$ where the minimum is attained at two or more distinct diagonal elements. To analyze this, let us express $f=\min \{f_1,\dots,f_n\}$, where $f_j$ is equal to $f$ except that the minimum over the diagonal elements in \eqref{eq:fqz} is replaced with the $j$-th diagonal element. A necessary condition for the possible non-differentiability of $f$ at $x_0$ is that there exists a subset $S \subseteq \{1,\dots,n\}$, depending on $x_0$ and such that (1) $\# S \geq 2$ (2) $f_i(x_0)=f(x_0)$ if $i \in S$ (3) $f_i(x_0) > f(x_0)$ if $i \not \in S$.

Suppose that $x_0$ is a local minimum of $f$ and that, in addition, we can find an associated set $S$ satisfying properties (1)-(2)-(3). Then, for all $x=(Q,Z)$ on the manifold and sufficiently close to $x_0$, if must hold $f(x_0) \leq f(x)$, which implies in particular $f_i(x_0) \leq f(x) \leq f_i(x)$ for all $i \in S$. Hence, $x_0$ is a local minimum of $f_i$ for all $i \in S$ (but not necessarily for $f_j$ when $j \not \in S$). Importantly, all $f_j$ are everywhere smooth functions on the real smooth manifold $U(n) \times U(n)$, because they depend polynomially on the real and imaginary parts of the entries of $Q$ and $Z$; hence, their real Riemannian gradient exists and vanishes at all their local minima. Denoting by $R_j(x)$ the Riemannian gradient of $f_j(x)$, we therefore have $R_i(x_0)=0$ for all $i \in S$. Fix a unit direction $v$ in the tangent space of $U(n) \times U(n)$ at $x_0$, and let $x(\epsilon)$ be the geodesics on $U(n) \times U(n)$ with $x(0)=x_0$, $\dot{x}(0)=v$. Then,
for all $i \in S$, $f_i(x(\epsilon))=  f_i (x_0)  + o(\epsilon)$ for sufficiently small $\epsilon$. Hence,
\[ \langle R(x_0),v \rangle = \lim_{\epsilon \rightarrow 0} \frac{f(x(\epsilon))-f(x_0)}{\epsilon} =  \lim_{\epsilon \rightarrow 0} \frac{\min_{i \in S}f_i(x(\epsilon))-f(x_0)}{\epsilon} = 0;\]
note that the limit is $0$ uniformly, as for sufficiently small $\epsilon$ it holds, for all $i$ and all $v$, $|f_i(x(\epsilon))-f(x_0)| < 2 \max_i \| H_i(x_0)\| \epsilon^2$ where $H_i$ is the Hessian of $f_i$.  We conclude that $f$ is real-differentiable at $x_0$ and that its differential is $0$.
    \end{proof}

 The (non-smooth) algorithm consisting of the minimization the function  \eqref{eq:newobjfun} via a Riemannian trust-region method \cite{ABG07} yields competitive solutions compared to other methods in the literature such as \cite{bhm,bora,glm}; moreover, we do not observe convergence issues in practice, and Theorem \ref{thm:lasttheorem} gives a clue that there is hope that in the future one could prove at least local convergence results; however, a full proof is beyond the scope of this paper. For readers who are nevertheless concerned about non-differentiability, a smoothed version of the algorithm is discussed in Section \ref{sec:smoothalgorithm}.}

% \begin{vn}
% I only made minor cosmetic changes to Remark 3.7, with one exception: I commented out the last sentence for brevity (it sounds unnecessary to repeat "you may want to go now to Section 6 if you want to know more on real pencils" -- I believe we already made it clear with the previous sentences that points to Lemma 6.1).

% Question: Is this the best place to write this Remark? Maybe it would be more natural right after Theorem 3.4
% \end{vn}

% \begin{fd}
% I agree with Vanni's rewriting. About moving the remark right after Theorem 3.4 also seems natural to me. Perhaps for the sake of readability I would move it right after the paragraph "This yields a method: we can minimize .... objective function." If not we should change "This yields..." by "Theorem 3.4 yields..." I will be happy with any option here.
% \end{fd}

\subsection{Riemannian optimization algorithm} \label{subsec:riemanprojec}

 The cost function $f(Q,Z)$ in \eqref{eq:minnewobjfun}, that we want to minimize over $U(n) \times U(n)$, can be written more explicitly as
\begin{equation}\label{eq:fqz}
      f(Q,Z) = \sum_{i>j} (|(QAZ)_{ij}|^2 + |(QBZ)_{ij}|^2) + \min_{1 \leq i \leq n} (|(QAZ)_{ii}|^2 + |(QBZ)_{ii}|^2).
\end{equation}
  
In order to express some results for $f(Q,Z)$ concisely, we define a function $L$ on any pencil of matrices $C+ \lambda D \in \C[\lambda]_1^{n \times n}$ as $L(C) + \lambda L(D) := (C-\mathcal{P}(C)) + \lambda ( D - \mathcal{P}(D))$, where $\mathcal{P}$ is as in Proposition~\ref{prop:nearest_tri}, i.e. $\mathcal{P}(C) + \lambda \mathcal{P}(D)$ is an upper triangular singular pencil nearest to $C+\lambda D$. In this setting, \eqref{eq:fqz} can be equivalently written as
\begin{equation} \label{eq:fqzL}
    f(Q,Z) = \| L(QAZ) \|_F^2 + \| L(QBZ) \|_F^2 . 
\end{equation}
Therefore, our algorithm aims to minimize this objective function over the Riemannian manifold $U(n) \times U(n)$. We use the MATLAB package Manopt \cite{manopt} version 7.1 as a toolbox to perform optimization on manifolds, and in particular Manopt's implementation of the trust-region method \cite{ABG07}. Necessary ingredients for the relevant Manopt's subroutines are the Riemannian gradient and the Riemannian Hessian \cite{rmanifolds} of the objective function. Next, we discuss how to obtain these quantities. \hide{If the Riemannian Hessian is not available, Manopt computes it from the gradient via a finite difference approach; the latter option has the advantage of not requiring the exact form of the Hessian, but numerically it may lead to suboptimal performance of the algorithm.} 

The Riemannian gradient of the objective function in the embedded manifold $U(n) \times U(n)  \subset \C^{n \times n} \times \C^{n \times n}$ can be computed by first computing the standard Euclidean gradient in the real vector ambient space $\C^{n \times n} \times \C^{n \times n} \cong \R^{4 n^2}$, and then performing an orthogonal projection onto the tangent space $T_{(Q,Z)} (U(n) \times U(n))$ of $U(n) \times U(n)$ at $(Q,Z)$ \cite[Proposition 3.61]{rmanifolds}. This orthogonal projection is defined with respect to the following real inner product in $\C^{n \times n} \times \C^{n \times n}$:
\begin{equation} \label{eq:innerprodunun}
\langle (A_1, A_2) \, , \, (B_1, B_2) \rangle := \mathrm{Re} \left(\mbox{trace} (A_1^* B_1) + \mbox{trace} (A_2^* B_2) \right),
\end{equation} 
for $ (A_1, A_2) ,  (B_1, B_2) \in \C^{n \times n} \times \C^{n \times n}$, which is induced by the real inner product $\langle A \, , \, B\rangle := \mathrm{Re} \left(\mbox{trace} (A^* B) \right)$ in $\C^{n \times n}$. The tangent space of a Cartesian product is the Cartesian product of the tangent spaces \cite[Proposition 3.20]{rmanifolds}, that is, 
\begin{align*}
    T_{(Q,Z)} (U(n) \times U(n)) = T_{Q} U(n) \times T_{Z} U(n).
\end{align*}
Moreover, the tangent space in $\C^{n \times n}$ of the manifold of unitary matrices at the unitary matrix $Q$ is $T_{Q} U(n)= \{ Q S: \ S = -S^*\}$ \cite[Lemma 3.2]{ant19}. It is a simple exercise to check that, with respect to the inner product $\langle A \, , \, B\rangle := \mathrm{Re} \left(\mbox{trace} (A^* B) \right)$, the orthogonal projection $\pi_{T_{Q} U(n)}$ of a matrix $M \in \C^{n \times n}$ onto the tangent space $T_{Q} U(n)$ is given by the skew-Hermitian projection $\text{skew}(M) = \frac{1}{2}(M - M^*)$ as $\pi_{T_{Q} U(n)} (M) = Q \ \text{skew}(Q^* M)$, and similarly for the tangent space $T_{Z} U(n)$. As a consequence, the orthogonal projection of $(A,B) \in \C^{n \times n} \times \C^{n \times n}$ onto the tangent space $T_{(Q,Z)} (U(n) \times U(n))$ is then $$\operatorname{Proj}_{(Q,Z)}(A,B) = (\pi_{T_{Q} U(n)}(A), \pi_{T_{Z} U(n)}(B)).$$
Thus, once we have computed the Euclidean gradient of the objective function (which we will do in Subsection \ref{subs:complexgrad}) we are fully equipped to then project it and obtain the Riemannian gradient.

In the case of a Riemannian submanifold $\mathcal{M}$ of a Euclidean space $\mathcal{E}$, the Riemannian Hessian of a smooth objective function $f : \mathcal{M} \rightarrow \R$ is \cite[Corollary~5.16]{rmanifolds}
\begin{align} \label{eq:rhess}
    \operatorname{Hess}_R f(x)[u]=\operatorname{Proj}_x(\mathrm{D} \bar{G}(x)[u]),
\end{align}
where $\operatorname{Proj}_x$ denotes the orthogonal projection onto $T_x \mathcal{M}$, the tangent space of $\mathcal{M}$ at $x$, and $\mathrm{D} \bar{G}(x)[u]$ denotes the directional derivative at $x \in \mathcal{M}$ of the vector field $\bar G$ in the direction $u \in T_x \mathcal{M}$. Here, $\bar G$ is any smooth vector field defined on a neighborhood of the manifold $\mathcal{M}$ in the embedding space $\mathcal{E}$ such that $\bar G(x) = \operatorname{Proj}_x(\nabla f(x))$ for all $x \in \mathcal{M}$, where $\nabla f(x)$ refers to the Euclidean gradient of the function $f$. In other words, $\bar G$ is a smooth extension of the Riemannian gradient of $f$. Our objective function $f(Q,Z)$ in \eqref{eq:fqz} is smooth almost everywhere, and so obeys this equality well in practice. 

The analytic expression of $\operatorname{Proj}_x(\nabla f(x))$ often provides a smooth extension of the Riemannian gradient in a natural manner. In practice however, computing $\operatorname{Proj}_x(\mathrm{D} \bar{G}(x)[u])$ via the expression of $\bar{G}(x)$ is not always the least laborious approach. In the case of the manifold $U(n) \times U(n)$, the Riemannian Hessian can be computed via { the Euclidean gradient and the Euclidean Hessian, $\operatorname{Hess}_E f(x)$, by Lemma~\ref{lemma:rhess}. An alternative expression for the second term in \eqref{eq.riehessia} can be obtained by following \cite[Section 4.3]{AMT13}.}

% \begin{vn}
%     I think that citing [2, Section 4.3] gives exactly zero useful information, but if you like it then it is certainly not worth arguing for.
% \end{vn}
% \begin{lauri}
%     Maybe the idea here is to please the referee? 
% \end{lauri}
In Lemma \ref{lemma:rhess} and throughout the paper, we define operations on ordered pairs as elementwise operations. In particular, we define the addition of two ordered pairs $(A_1,A_2), (B_1,B_2) \in \C^{n \times n} \times \C^{n \times n}$ by $(A_1,A_2) + (B_1,B_2) := (A_1 + B_1,A_2 + B_2)$, their multiplication by $(A_1,A_2) \cdot (B_1,B_2) := (A_1 B_1,A_2 B_2)$, and the operation of taking the conjugate transpose by $(A_1,A_2)^* := (A_1^*,A_2^*)$. Moreover, the orthogonal projections are considered with respect to the inner product \eqref{eq:innerprodunun}.

\begin{lemma} \label{lemma:rhess}
Let $f: U(n) \times U(n) \longrightarrow \R$ be a smooth function. Then, for any $x \in U(n) \times U(n)$ and $u \in T_x \,( U(n) \times U(n) )$, it holds that
\begin{align} \label{eq.riehessia}
    \operatorname{Hess}_R f(x)[u] = \operatorname{Proj}_x \left(\operatorname{Hess}_E f(x)[u] - \frac{1}{2} u \, (x^* \nabla f(x) + \nabla f(x)^* x) \right),
\end{align}
where the Euclidean gradient $\nabla f(x) \in \C^{n \times n} \times \C^{n \times n}$ at $x = (x_1,x_2)$ is the ordered pair defined as 
$$
\nabla f(x) = ( \nabla_{x_1} f(x_1,x_2) \, , \, \nabla_{x_2} f(x_1,x_2))
$$
with $\nabla_{x_i} f(x_1,x_2)$ the Euclidean gradient with respect to the argument $x_i \in \C^{n \times n}$ for $i=1,2$.
\end{lemma}
\begin{proof}
Let us use $\operatorname{Proj}_x(\nabla f(x))$ as the choice for the smooth extension $\bar G (x)$ in \eqref{eq:rhess}. Substituting this into (\ref{eq:rhess}) yields
\begin{align*}
    \operatorname{Hess}_R f(x)[u] &= \operatorname{Proj}_x(\mathrm{D} \bar{G}(x)[u]) = \operatorname{Proj}_x(\mathrm{D} \operatorname{Proj}_x(\nabla f(x))[u]) \\
    &= \operatorname{Proj}_x(\mathrm{D} \nabla f(x)[u]) - \operatorname{Proj}_x(\mathrm{D} (I-\operatorname{Proj}_x)(\nabla f(x))[u]).
\end{align*}
Since $\operatorname{Hess}_E f(x)[u] = \mathrm{D} \nabla f(x)[u]$, it only remains to show that 
\begin{align} \label{eq:proj_lemma}
    \operatorname{Proj}_x(\mathrm{D} (I-\operatorname{Proj}_x)(\nabla f(x))[u]) = \operatorname{Proj}_x(\frac{1}{2} u (x^* \nabla f(x) + \nabla f(x)^*x)).
\end{align}
Let us recall that
\begin{align*}
    (I-\operatorname{Proj}_x)(\nabla f(x)) = \frac{1}{2} x \, ( x^* \nabla f(x) + \nabla f(x)^* x).
\end{align*}
Then,
\begin{align*}
    \mathrm{D} (I-\operatorname{Proj}_x)(\nabla f(x))[u] &= \frac{1}{2} x ( x^* \operatorname{Hess}_E f(x)[u] + \operatorname{Hess}_E f(x)[u]^* x) \\ &\;\;\; +  \frac{1}{2} ( u \nabla f(x)^* x + x \nabla f(x)^* u)\\
    &= (I-\operatorname{Proj}_x)(\operatorname{Hess}_E f(x)[u]) + \frac{1}{2} ( u \nabla f(x)^* x + x \nabla f(x)^* u).
\end{align*}
When we apply $\operatorname{Proj}_x$ again, the first term becomes zero. Hence, the left hand side of (\ref{eq:proj_lemma}) becomes
\begin{align*}
     \operatorname{Proj}_x(\mathrm{D} (I-\operatorname{Proj}_x)(\nabla f(x))[u]) = \operatorname{Proj}_x \left(\frac{1}{2} ( u \nabla f(x)^* x + x \nabla f(x)^* u)\right) \\
     = \frac{1}{4} (u \nabla f(x)^* x - \nabla f(x) u^* x) + \frac{1}{4}(x \nabla f(x)^* u - x u^* \nabla f(x)) .
\end{align*}
On the other hand the right hand side of (\ref{eq:proj_lemma}) is
\begin{align*}
    \operatorname{Proj}_x \left(\frac{1}{2} u (x^* \nabla f(x) + \nabla f(x)^*x)\right) & = \frac{1}{4} ( u x^* \nabla f(x) - x \nabla f(x)^* x u^* x) \\ & \;\;\; + \frac{1}{4} ( u \nabla f(x)^* x - \nabla f(x) u^* x).
\end{align*}
As $u \in T_{x} (U(n) \times U(n)) = \{ x_1 S: \ S = -S^*\} \times  \{ x_2 S: \ S = -S^*\}$, where $x = (x_1, x_2)$, it holds that $x u^* = - u x ^*$. Then, $u x^* \nabla f(x) = - xu^* \nabla f(x)$ and $x \nabla f(x)^* x u^* x = - x \nabla f(x)^*u$. Thus, the equality in (\ref{eq:proj_lemma}) holds.
\end{proof}

Lemma~\ref{lemma:rhess} lets us compute the Riemannian Hessian via the Euclidean Hessian and the Euclidean gradient. In order to use this approach, we derive an expression for the Euclidean Hessian of the objective function in Subsection~\ref{subs:complexhess}.

\subsection{The Euclidean gradient of the objective function}\label{subs:complexgrad}

In this subsection we aim to compute the Euclidean gradient 
\begin{align*} 
    \nabla f(Q,Z) := (\nabla_{Q} f(Q,Z), \nabla_{Z} f(Q,Z)) \in \C^{n \times n} \times \C^{n \times n},
\end{align*} 
where the right hand side is an ordered pair consisting of the gradients with respect to the arguments $Q$ and $Z$. Note that, generally, $f(Q,Z)$ is not differentiable; however, it is almost everywhere differentiable, so that in practice we can expect the gradient to exist when evaluating it in our numerical algorithm.

 To this goal, motivated by \eqref{eq:fqzL}, let us first consider a univariate function $g(Q) := \| L(QAZ) \|_F^2$ which can in turn  be seen as a composition $g = g_1 \circ g_2$, where $g_1(X) = \| X \|_F^2$ and $g_2(X) = L(XAZ)${, where $L$ is piecewise linear}. The directional derivatives of $g_1$ and $g_2$ {(whenever $XAZ$ is in the interior of a linear piece of $L$)} are, respectively, 
\begin{align*}
    \mathrm{D} g_1(Q)[X] &= 2 \, \mathrm{Re} (\tr (Q^*X)), \\
    \mathrm{D} g_2(Q)[X] &= L(XAZ).
\end{align*}
By the chain rule, it follows that
\begin{align*}
    \mathrm{D} g(Q)[X] &= (\mathrm{D} g_1 (g_2(Q)) \circ \mathrm{D} g_2 (Q))[X]\\
    &= 2 \, \mathrm{Re} (\tr (L(QAZ)^*L(XAZ))) \\
    &= 2 \, \mathrm{Re} ( \tr ((L(QAZ)(AZ)^*)^* X )). 
    % &= \cvec (2 L(QAZ) (AZ)^*)^T \cvec (X).
\end{align*}

The gradient $\nabla_{Q} g(Q)$ is the unique vector which satisfies $$\mathrm{D} g(Q)[X] = \langle\nabla_{Q} g(Q),X \rangle = \mathrm{Re} (\tr (\nabla_{Q} g(Q)^* X)) $$ for all $X$. Hence, the gradient is $\nabla_{Q} g(Q) = 2 L(QAZ) (AZ)^*$.

An analogous result holds for $\| L(QBZ) \|_F^2$ and so
\begin{align} \label{eq:gradient}
    \nabla_{Q} f(Q,Z) &= 2 L(QAZ) (AZ)^* + 2 L(QBZ) (BZ)^*.
\end{align}

We can compute $\nabla_{Z} f(Q,Z)$ from the formula for $\nabla_{Q} f(Q,Z)$ by transposing the matrices in an appropriate manner, which yields
\begin{align*}
    \nabla_{Z} f(Q,Z) &= 2 (QA)^* L(QAZ) + 2 (QB)^* L(QBZ).
\end{align*}

\subsection{The Euclidean Hessian of the objective function} \label{subs:complexhess}
The Euclidean Hessian is defined via the directional derivative of the gradient \cite[p. 23]{rmanifolds}. In other words, we are interested in
\begin{align} \label{eq:euclideanhess}
    \operatorname{Hess}_E f (Q,Z)[d] := \mathrm{D} \nabla f (Q,Z)[d] = (\mathrm{D} \nabla_Q f (Q,Z)[d],\mathrm{D} \nabla_Z f (Q,Z)[d]),
\end{align}
where $d$ is the vector representing the direction. Here, $d$ is an ordered pair $d := (d_Q, d_Z) \in \C^{n \times n} \times \C^{n \times n}$, where $d_Q$ and $d_Z$ are the matrix directions for the matrices $Q$ and $Z$, respectively. The rightmost expression is an ordered pair consisting of the directional derivatives of $\nabla_{Q} f(Q,Z)$ and $\nabla_{Z} f(Q,Z)$, respectively.

Consider again the function $g(Q) := \| L(QAZ) \|_F^2$. Its gradient was computed to be $\nabla_{Q} g(Q,Z) = 2 L(QAZ) (AZ)^*$ in Subsection \ref{subs:complexgrad}. The directional derivative with respect to $Q$ of $\nabla_{Q} g(Q,Z)$ is 
\begin{align*}
\mathrm{D}^Q \nabla_{Q} g(Q,Z) [X] = 2 L(XAZ)(AZ)^*.  
\end{align*}
Moreover, the directional derivative with respect to $Z$ is 
\begin{align*}
\mathrm{D}^Z \nabla_{Q} g(Q,Z) [X] = 2 L(QAX)(AZ)^* + 2 L(QAZ)(AX)^*.    
\end{align*}
The directional derivative of $\nabla_{Q} g(Q,Z)$ is the sum of these two evaluated at the matrix directions  $d_Q$ and $d_Z$, respectively:
\begin{align*}
    \mathrm{D} \nabla_{Q} g (Q,Z)[d_Q,d_Z] = 2 L(d_Q AZ)(AZ)^* + 2 L(QA d_Z)(AZ)^* + 2 L(QAZ)(A d_Z)^*.
\end{align*}
An analogous result holds for the term $\| L(QBZ) \|_F^2$ in \eqref{eq:fqzL}, and so the directional derivative of  $\nabla_{Q} f(Q,Z)$ is
\begin{align*}
    \mathrm{D} \nabla_{Q} f (Q,Z)[d_Q,d_Z] &= 2 L(d_Q AZ)(AZ)^* + 2 L(QA d_Z)(AZ)^* + 2 L(QAZ)(A d_Z)^*\\
    &+ 2 L(d_Q BZ)(BZ)^* + 2 L(QB d_Z)(BZ)^* + 2 L(QBZ)(B d_Z)^*.
\end{align*}
The directional derivative for $\nabla_{Z} f(Q,Z)$ can be computed in a similar way, which yields
\begin{align*}
    \mathrm{D} \nabla_{Z} f(Q,Z)[d_Q,d_Z] &= 2 (QA)^*L(d_Q AZ) + 2 (d_Q A)^* L(QAZ) + 2 (QA)^* L(QA d_Z) \\
    &+ 2 (QB)^*L(d_Q BZ) + 2 (d_Q B)^* L(QBZ) + 2 (QB)^* L(QB d_Z).
\end{align*}

{ 
\section{Smoothed algorithm}\label{sec:smoothalgorithm} As discussed in Subsection \ref{subsec:reformulating}, the objective function given in \eqref{eq:fqz} is non-smooth. % due to {\colb the addend including} the \texttt{min}-function. 
This is potentially problematic because traditional smooth optimization techniques may not always converge to stationary points on non-smooth problems, as discussed in \cite{asl-overto-2020} and the references therein. %Thus, while we anticipate convergence for almost every choice of $A,B$ in Problem \ref{problem}, it is possible that inputs exist for which the algorithm developed in Section \ref{sec:complex} does not converge to a stationary point.

To address this concern, we can smoothen the objective function by replacing the \texttt{min}-function with the Boltzmann operator, which is one of the most commonly used softmin operators \cite{asadi17}. This operator depends on a real parameter $\alpha$ and for $n$ real variables $x_1, \ldots, x_n$ is defined as
$$
\mathcal{S}_\alpha\left(x_1, \ldots, x_n\right)=\frac{\sum_{i=1}^n x_i e^{\alpha x_i}}{\sum_{i=1}^n e^{\alpha x_i}}.
$$
The Boltzmann operator $\mathcal{S}_\alpha \left(x_1, \ldots, x_n\right)$ is infinitely differentiable as a function of $\left(x_1, \ldots, x_n\right)$ for any value of the real parameter $\alpha$. The key observation in our context is that $\mathcal{S}_\alpha \left(x_1, \ldots, x_n\right) \rightarrow \min \left(x_1, \ldots, x_n\right)$ as $\alpha \rightarrow-\infty$. Hence, a smooth version of the objective function \eqref{eq:fqz} can be given as 
\begin{equation}\label{eq:smooth}
 f_\alpha  (Q,Z) = \sum_{i>j} (|(QAZ)_{ij}|^2 + |(QBZ)_{ij}|^2) + \mathcal{S}_\alpha\left(x_1, \ldots, x_n\right),
\end{equation}
where $x_i = |(QAZ)_{ii}|^2 + |(QBZ)_{ii}|^2$ and $\alpha$ is a negative value of a suitable magnitude. The Riemannian gradient and Hessian for the smooth objective function \eqref{eq:smooth} can be derived following the approach in Section \ref{sec:complex}. That is, computing first the Euclidean gradient and then orthogonally projecting onto the tangent space $T_{(Q,Z)} (U(n) \times U(n))$, and computing the Euclidean Hessian and then using Lemma \ref{lemma:rhess}. The computations of the Euclidean gradient and Hessian of $f_\alpha (Q,Z)$ are considerably more involved than for $f(Q,Z)$ in \eqref{eq:fqz}. Thus, we omit the details and limit ourselves to present their mathematical expressions in Appendix \ref{sec:gradient_smooth}.

It should be borne in mind that, once a minimizer $(Q_0,Z_0)$ of $f_\alpha (Q,Z)$ is computed, the distance to singularity should be computed by evaluating the square root of $f(Q_0,Z_0)$, with $f(Q,Z)$ as in \eqref{eq:newobjfun}, and that the computed  nearest singular pencil is obtained as in the second part of Theorem \ref{thm:obj_function}.

 We emphasize that for each fixed value of $\alpha$, the function $f_\alpha (Q,Z)$ is everywhere smooth on the real manifold $U(n) \times U(n)$, because $|(QAZ)_{ij}|^2 + |(QBZ)_{ij}|^2$ is a polynomial function of the real and imaginary parts of the entries of $Q$ and $Z$  for each $(i,j)$. Moreover, the manifold $U(n) \times U(n)$ is compact. Thus, Corollary 7.4.6 and Theorem 7.4.4 in \cite{AMS08} guarantee that the Riemannian trust-region method in Manopt \cite{manopt} converges globally to a stationary point of $f_\alpha (Q,Z)$ \cite[Corollary 4.6]{ABG07}, which in computational practice is a local minimizer except for very special situations  \cite[Section 7.4.3]{AMS08}.

In the numerical experiments presented in Section \ref{sec:numerical}, the standard non-smoothed algorithm from Section \ref{sec:complex} always converged to a stationary point for a large volume of randomly generated inputs, as well as on benchmark problems previously appeared in the literature. This behaviour is coherent with the expectations, as the objective function is  differentiable almost everywhere and always in its local minima. In Subsection \ref{sec:smooth_numerical}, we will see that the smoothed algorithm presented in this section can be significantly slower than the non-smoothed version. For these reasons, we opted for using the non-smoothed version as the primary algorithm for the numerical experiments and we recommend to use the non-smooth algorithm in practice. Nevertheless, for users concerned with the lack of provable convergence, the smoothed version is a viable alternative since it is not dramatically slower.
}

{
\subsection{A third algorithm}\label{sec:third}
Section \ref{sec:fixminindex} presents an algorithm for a variant of the nearest singular pencil problem, namely, the problem of finding the nearest singular pencil with a prescribed minimal index. It turns out that the algorithm of Section \ref{sec:fixminindex} can also be used to construct an alternative approach in providing a smooth algorithm for the general nearest singular pencil problem, by separately optimizing over all the posssible branches of the minimum function.

More in detail, observe that a global minimizer $(Q_0,Z_0)$ of $f(Q,Z)$ in \eqref{eq:newobjfun} or \eqref{eq:fqz} satisfies $\min_{1 \leq i \leq n} (|(Q_0AZ_0)_{ii}|^2 + |(Q_0BZ_0)_{ii}|^2)=|(Q_0AZ_0)_{kk}|^2 + |(Q_0BZ_0)_{kk}|^2$ for some index $k \in \{1, \ldots, n\}$. Thus, one can find it by instead minimizing over $U(n) \times U(n)$ the everywhere smooth function
\begin{equation} \label{eq:fkfirst}
f_k (Q,Z) = \sum_{i>j} (|(QAZ)_{ij}|^2 + |(QBZ)_{ij}|^2) + (|(QAZ)_{kk}|^2 + |(QBZ)_{kk}|^2).
\end{equation}
While the precise index $k$ is unknown in advance, one can still minimize $f(Q,Z)$ by (a) separately computing $n$ minimizers of the smooth functions $f_1 (Q,Z), \ldots, f_n (Q,Z)$ over $U(n) \times U(n)$, and (b) picking the one associated with the smallest minimum. Subsection \ref{sec:minindex_numerical} demonstrates that this approach performs equivalently (in a sense) to running the non-smooth algorithm $n$ times. Additionally, each minimization problem $\min f_i (Q,Z)$ solves the problem of finding a nearest singular pencil with a specified (right) minimal index equal to $i-1$; the proof of this fact requires several new theoretical results, which we develop in Section \ref{sec:fixminindex}. The same argument used above for $f_\alpha (Q,Z)$ in \eqref{eq:smooth} proves that the Riemannian trust-region method in Manopt converges globally to a stationary point of $f_i(Q, Z)$ for each $i =1, \ldots , n$, which in practice is a local minimizer except in very special cases.
}

%{\colb 
\section{Nearest singular pencil with a specified minimal index}\label{sec:fixminindex}
% A minimizer $(Q_0,Z_0)$ of $f(Q,Z)$ in \eqref{eq:newobjfun} (or, in \eqref{eq:fqz}) obviously satisfies $|(Q_0AZ_0)_{kk}|^2 + |(Q_0BZ_0)_{kk}|^2 = \min_{1 \leq i \leq n} (|(Q_0AZ_0)_{ii}|^2 + |(Q_0BZ_0)_{ii}|^2)$ for some index $k \in \{1, \ldots, n\}$. Thus, this minimizer would be found if instead of minimizing $f(Q,Z)$ one minimizes over $U(n) \times U(n)$ the function 
% \begin{equation} \label{eq:fkfirst}
% f_k (Q,Z) =  \sum_{i>j} (|(QAZ)_{ij}|^2 + |(QBZ)_{ij}|^2) + (|(QAZ)_{kk}|^2 + |(QBZ)_{kk}|^2),
% \end{equation}
% which is a polynomial in the real and imaginary parts of the entries of $Q$ and $Z$ and, thus, smooth. Of course, the precise index $k$ is not known in advance but if one solves separately over $U(n) \times U(n)$ the $n$ smooth minimization problems $\min f_1 (Q,Z), \ldots$, $\min f_n (Q,Z)$ and gets as minimizers $(Q_1,Z_1), \ldots , (Q_n,Z_n)$, respectively, then a minimizer $(Q_0,Z_0)$ of $f(Q,Z)$ can be obtained by $f(Q_0,Z_0) = \min_{1 \leq i \leq n} f_i(Q_i,Z_i)$. Subsection \ref{sec:minindex_numerical} demonstrates that the performance of this approach is equivalent (in a particular sense) to the performance of the non-smoothed algorithm. Moreover, we prove in this section that each minimization problem $\min f_i (Q,Z)$ solves the problem of finding a nearest singular pencil with a specified (right) minimal index equal to $i-1$ {\colb and that this problem is well posed}. This requires to first develop several new theoretical results.}

{Recall that the normal rank of a matrix pencil $A + \lambda B \in \C[\lambda]_1^{n \times n}$ refers to its rank over the field of rational functions $\C(\lambda)$. If the pencil has normal rank $n - 1$, there exists a one-dimensional space of rational vectors in its right null space. Some of these vectors are polynomial, and the smallest degree among these nonzero polynomial vectors is called the right minimal index. See \cite{forney} for a more general definition of minimal indices.}

Let $\mathcal{S}_n^k \subset \C[\lambda]_1^{n \times n}$ denote the set of singular pencils with normal rank $n-1$ and right minimal index $k$. The set $\mathcal{S}_n^k$ is not closed and consequently a minimizer of the distance function \eqref{eq:pencildist} might not be well defined. Instead, we can consider the distance to the set $\mathcal{S}_n^k$ (i.e., the infimum of the distance function), and construct a pencil in $\mathcal{S}_n^k$ that is arbitrarily close to this distance. This is stated formally in Problem~\ref{problem_minimal_index}.

\begin{problem}\label{problem_minimal_index}
Given a pencil $A+\lambda B \in \C[\lambda]_1^{n \times n}$, a nonnegative integer $k \leq n-1$, and a parameter $\epsilon > 0$, compute the value $$\inf_{S+\lambda T \in \mathcal{S}_n^k} d(A+\lambda B,S+\lambda T),$$ and construct a pencil $S_\epsilon+\lambda T_\epsilon \in \mathcal{S}_n^k$ such that

$$d(A+\lambda B,S_\epsilon+\lambda T_\epsilon) \leq \inf_{S+\lambda T \in \mathcal{S}_n^k} d(A+\lambda B,S+\lambda T) + \epsilon.$$

\end{problem}

We approach Problem~\ref{problem_minimal_index} by finding a minimizer over the closure of $\mathcal{S}_n^k$, denoted by $\overline{\mathcal{S}_n^k}$, and then finding an arbitrarily small perturbation that brings this minimizer back to $\mathcal{S}_n^k$ (note that $\overline{\mathcal{S}_n^k}$ being closed guarantees the existence of a minimizer).

Let $\text{Schur}_n^{k} \subset \C[\lambda]_1^{n \times n}$ denote the set of pencils for which there exists a generalized Schur form in which the $k$th diagonal element is zero. Later in this section, we prove in Theorem~\ref{thm:arbitrarily_close} a result that yields $\text{Schur}_n^{k+1} \subseteq \overline{\mathcal{S}_n^k}$, and we obtain in Lemma \ref{lemm:auxxk1} a result that in turn implies $\mathcal{S}_n^k \subseteq \text{Schur}_n^{k+1}$. Finally, we show with Theorem~\ref{thm:closure} that $\text{Schur}_n^{k+1}$ is closed. Collectively, these three facts imply that $\text{Schur}_n^{k+1} = \overline{\mathcal{S}_n^k}$ (Corollary \ref{cor:si}). We then solve Problem~\ref{problem_minimal_index} by computing 
\begin{align}\label{eq:riproblem_k}
  \min_{S+\lambda T \in \text{Schur}_n^{k+1}} d(A+\lambda B,S+\lambda T),  
\end{align} 
and use the proof of Theorem~\ref{thm:arbitrarily_close} to find an arbitrarily small perturbation that brings the minimizer to $\mathcal{S}_n^k$. It is worth observing that the set of $n\times n$ singular complex pencils is precisely $\bigcup_{0\leq k \leq n-1} \overline{\mathcal{S}_n^k}$  \cite[Theorem 3.2]{deterandopico08}. { Observe that solving \eqref{eq:riproblem_k} for $k=0, \ldots, n-1$ is precisely the strategy for solving Problem \ref{problem} discussed in Subsection \ref{sec:third}.}

%This will be illustrated numerically in Subsection %\ref{sec:minindex_numerical}.

In the following, we follow the convention that $L(\infty)$, where $L$ is a polynomial matrix over $\C$ having grade $g$, denotes the $g$-th degree coefficient of $L$, which is allowed to be zero. For brevity, within this section we often omit the dependence on the variable in polynomial matrices and pencils. Recall that a singular polynomial matrix $L$ has, { a nontrivial} rational right null space, that is denoted by $\ker L$ and  that $\ker L$ has minimal bases \cite{forney}. These minimal bases consist of polynomial vectors whose degrees are the right minimal indices of $L$. Moreover, $L(\lambda)$ stands for the evaluation of $L$ at $\lambda \in \mathbb{C}$. Given a constant matrix $M$, $\nu(M)$ denotes the nullity of $M$ i.e. the dimension of the right null space of $M$. Finally, $\overline{\C}:=\C \cup \{ \infty \}$.

\begin{lemma}\label{lem1}
Let $T$ be a singular $n \times n$ upper triangular pencil such that (1) the normal rank of $T$ is $n-1$ and (2) $T$ has no eigenvalues. Let the uppermost zero diagonal element of $T$ be located at position $(k+1,k+1)$. Then, the right minimal index of $T$ is $k$.
\end{lemma}

\begin{proof}
$T$ has the form
\begin{equation}\label{eq:niceform}
T = \begin{bmatrix}
T_1 & R_1 & R_2\\
0 & 0 & R_3\\
0 & 0 & T_2
\end{bmatrix},
\end{equation}
where $T_1$ is a $k \times k$ triangular invertible pencil, $T_2$ is an $(n-k-1) \times (n-k-1)$ triangular pencil, and sizes of the other (possibly) full blocks are as follows: $R_1$ is $k \times 1$, $R_2$ is $k \times (n-k-1)$ and $R_3$ is $1 \times (n-k-1)$. Note that $\begin{bmatrix}
R_3 \\
T_2
\end{bmatrix}$ must have full column rank, otherwise $\dim \ker T \geq 2$ contradicting assumption (1). Then, every rational vector in $\ker T$ must have the form $\begin{bmatrix}
p\\
q\\
0
\end{bmatrix}$ where $q$ is $1\times 1$ and $p$ is $k \times 1$. Moreover, it must hold that $T_1 p = - R_1 q$. Define $A_1=\adj(T_1)$, then a polynomial solution is $q=\det (T_1)$, $p = -A_1 R_1$. Moreover, we henceforth define the grades of $A_1$ and $A_1 R_1$ to be $k-1$ and $k$ respectively.  To show that the polynomial vector $\begin{bmatrix}
-A_1R_1\\
\det(T_1)\\
0
\end{bmatrix}$
is a minimal basis of $\ker T$ we use \cite[p. 495, Main Theorem, item 2]{forney}. For this purpose, let $(\mu,v,w)$ be any eigentriple of $T_1$. We claim that $\mu$ has geometric multiplicity $1$.  In this case, $A_1(\mu)=-c v w^*$ for some $0 \neq c \in \C$ \cite[Section 0.8.2]{hornjohnson}. Hence, $-A_1(\mu)R_1(\mu)=c v w^* R_1(\mu)$. (Recall that, if $\mu=\infty$, $p(\infty)$ is the degree $k$ coefficient of $p$ and $A_1(\infty)$ is the degree $k-1$ coefficient of $A_1$.) Suppose $A_1(\mu)R_1(\mu)=0$, then $R_1(\mu) \in \spann(w)^\perp = \mathrm{colspace}(T_1(\mu))$, which implies \[\exists a : R_1(\mu)=T_1(\mu) a \Rightarrow \nu(\begin{bmatrix}
T_1(\mu) & R_1(\mu)
\end{bmatrix}) \geq 2\]
(to see the last step, note that $\begin{bmatrix}
v\\
0
\end{bmatrix}$ and $\begin{bmatrix}
a\\
-1
\end{bmatrix}$ must be linearly independent)
which is again a contradiction as it leads to the conclusion that $\mu$ is an eigenvalue of $T$. Therefore, $A_1(\mu)R_1(\mu) \ne 0$ and $\begin{bmatrix}
-A_1R_1\\
\det(T_1)\\
0
\end{bmatrix}$ is a minimal basis of $\ker T$ with degree $k$.

To conclude the proof, it remains to prove the claim that $\mu$ has geometric multiplicity $1$. To this goal, let $\mu \in \overline{\C}$ be an eigenvalue of $T_1$: then its geometric multiplicity must be $1$, since otherwise
\begin{equation}\label{eq:ohno}
    \nu(T(\mu))  \geq \nu(T_1(\mu)) \geq 2.
\end{equation} 
Indeed, if $B$ is a full rank matrix such that $T_1(\mu) B=0$ then $T(\mu) \begin{bmatrix}
B\\ 0 \\
0
\end{bmatrix}=0$. Hence, \eqref{eq:ohno} cannot be true as it would imply that $\mu$ is an eigenvalue of $T$ contradicting assumption (2).
\end{proof}

\begin{theorem} \label{thm:arbitrarily_close}
Let $U$ be a singular $n \times n$ upper triangular pencil with a zero diagonal element in the $(k+1,k+1)$ entry. Then, there exists a singular upper triangular pencil $T$ satisfying the following properties: (i) $T$ has normal rank $n-1$ (ii) $T$ has no eigenvalues (iii) $T$ has exactly one zero diagonal element located in the $(k+1,k+1)$ entry (iv) The unique right minimal index of $T$ is equal to $k$ (v) $T$ is arbitrarily close to $U$.
\end{theorem}

\begin{proof}
Let $\epsilon > 0$. We first perturb all the diagonal elements in $U$, other than in the $(k+1,k+1)$ position, so that they are all nonzero and all evaluate to zero at a distinct element of $\overline{\C}$. This can be achieved by perturbing each element by a perturbation of norm $\leq \frac{\epsilon}{3 \sqrt{n-1}}$ and creates a pencil $P$ such that $\|P-U\|_F \leq \frac{\epsilon}{3}$. 

Obviously $P$ has normal rank $n-1$, but it might have eigenvalues, and we still need to avoid the latter event. To this goal, note that $P$ has the form \eqref{eq:niceform}, and let us call the corresponding blocks $T_1,T_2,R_1,R_2,R_3$ as in \eqref{eq:niceform}. Moreover, observe that $T_1,T_2$ are both regular. For every $\mu \in \overline{\C}$, $\rank P(\mu) \geq \rank T_1(\mu) + \rank T_2(\mu)$, and hence an eigenvalue of $P$ must also be an eigenvalue of either $T_1$ or $T_2$ (and not both since the diagonal elements of the $T_i$ all evaluate to $0$ at distinct points in the Riemann sphere). Denote by $\mu_i$, $i=1,\dots,k$, the eigenvalues of $T_1$. Since $T_2(\mu_i)$ is invertible, it holds that $\rank P(\mu_i) = \rank \begin{bmatrix}
T_1(\mu_i) & R_1(\mu_i)
\end{bmatrix} + n-k-1$. Hence, $\mu_i$ is an eigenvalue of $P$ if and only if it is an eigenvalue of $\begin{bmatrix}
T_1 & R_1
\end{bmatrix}$. We can now make sure that this does not happen by only perturbing the constant coefficient of $R_1$. Write $R_1 = a + b x$ for some constant vectors $a,b \in \C^k$. For each $i$, let $\mathcal{A}_i$ denote the $(k-1)$-dimensional affine space $\mathcal{A}_i:=- \mu_i b + \mathrm{colspace}\, (T_1 (\mu_i))$.  For a fixed $i$, the condition $a \in \mathcal{A}_i$ is described by an affine equation in the entries of $a$, say $\ell_i(a)=0$. Hence, the polynomial equation $\prod_{i=1}^k \ell_i(a) = 0$ describes the condition that $a \in \mathcal{A}_i$ for at least one value of $i$. This is a proper Zariski closed set, hence its complement is non-empty and Zariski open, which means that it is open and dense in the Euclidean topology. Thus, we can perturb $a$, and therefore $P$, by a perturbation of norm bounded by $\frac{\epsilon}{3}$ to guarantee that there are no eigenvalues in $\begin{bmatrix}
T_1 & R_1
\end{bmatrix}$. Next, we repeat the same procedure by perturbing $R_3$, and hence $P$, again by a perturbation of norm $\leq \frac{\epsilon}{3}$ to make sure that $\begin{bmatrix}
    R_3\\
    T_2
\end{bmatrix}$ has no eigenvalues. We have thus constructed an upper triangular pencil $T$ with properties (i), (ii), (iii); moreover, property (v) follows from
\[ \| T - U \|_F \leq \|T-P\|_F + \|P-U\|_F \leq \frac{2 \epsilon}{3} + \frac{\epsilon}{3} = \epsilon.\]
Finally, we apply Lemma \ref{lem1} to $T$ to show that it has minimal index $k$, proving property (iv).
\end{proof}

The generalized Schur form of a pencil is not unique. However, for pencils having one of its right minimal indices equal to $k$ there always exist some generalized Schur forms with their $(k+1 , k+1)$ entry equal to zero. This is proved in Lemma \ref{lemm:auxxk1} below.

\begin{lemma} \label{lemm:auxxk1} Let $W$ be a singular $n\times n$ pencil having one of its right minimal indices equal to $k$, $0 \leq k \le n-1$. Then, there exist $Q, Z \in U(n)$ such that $Q W Z$ is upper triangular and its $(k+1, k+1)$ entry is equal to zero.
\end{lemma}
\begin{proof}
The right minimal index $k$ of $W$ implies that the complex Kronecker canonical form of $W$  (which is unique up to permutation of its blocks, see \cite{gant59}) has one block equal to
$$
K_k = \begin{bmatrix}
                \lambda & 1 &  &  \\
                        & \ddots & \ddots &  \\
                 &  & \lambda & 1
              \end{bmatrix} \in \C[\lambda]_1^{k \times (k+1)}.
$$
Hence, there exist invertible constant matrices $G, H \in \C^{n \times n}$ such that $W = G (K_k \oplus \widetilde{W} ) H$, where $K_k \oplus \widetilde{W}$ is a Kronecker canonical form of $W$ (note that if $k=0$, then $K_0$ simply appends a zero column on the left of $\widetilde{W}$). Observe that it is possible to choose the order of the blocks in the Kronecker canonical form such that the $n\times n$ pencil $K_k \oplus \widetilde{W}$ is upper triangular and that its $(k+1, k+1)$ entry is equal to zero. 
Let us now consider $QR$ factorizations $G = Q^* R_G$ and $H = R_H Z^*$ with $Q,Z \in U(n)$ and $R_G, R_H$ upper triangular. Thus
$$
W = Q^* \underbrace{\Big( R_G (K_k \oplus \widetilde{W} )R_H \Big) }_{:=T} Z^*,
$$
and the pencil $T=QWZ$ is upper triangular and has its $(k+1, k+1)$ entry equal to zero.
\end{proof}

We next show that $\text{Schur}_n^{k}$ is a closed set.

\begin{theorem} \label{thm:closure}
Let $1 \leq k \leq n$. The set $\text{Schur}_n^{k}$ is closed.
\end{theorem}
\begin{proof}
    We show that $\text{Schur}_n^{k}$ contains all of its accumulation points. Let $U \in \C[\lambda]_1^{n \times n}$ be an accumulation point of $\text{Schur}_n^{k}$. Then, there exists a sequence $(V_i)_{i\in \N},$ $V_i \in \text{Schur}_n^{k}$ such that $\lim_{i \rightarrow \infty} \| V_i - U \|_F = 0$.  For each $V_i$ there exists a triple $(Q_i,T_i,Z_i)$ such that $Q_i,Z_i$ are unitary and $Q_i^* V_i Z_i^* = T_i$ is upper triangular with a zero in position $(k,k)$. 
    
    Let $\text{Upp}^k_n$ denote the set of upper triangular pencils of size $n\times n$ with a zero in position $(k,k)$. It is straightforward to see that $\text{Upp}^k_n$ is a Zariski closed set. For any $\epsilon > 0$, the intersection of the closed ball $B_{\|U\|_F + \epsilon}(0)$ with $\text{Upp}^k_n$ is compact. Moreover, there exists an index $j$ such that $(T_i)_{i>j}$ is a sequence in $B_{\|U\|_F + \epsilon}(0) \cap \text{Upp}^k_n$ as $|\| T_i \|_F - \| U \|_F| = |\| V_i \|_F - \| U \|_F| \leq \| V_i - U \|_F$ and $\lim_{i \rightarrow \infty} \| V_i - U \|_F = 0$. As the set of unitary matrices is also compact (it is Zariski closed and bounded), and the Cartesian product of compact sets is compact, then $(Q_i,T_i,Z_i)_{i>j}$ is a sequence in a compact set, and hence has a converging subsequence whose limit belongs to the set itself. Let $(Q_{i_l},T_{i_l},Z_{i_l})_{l\in \N}$ denote such a sequence, where $j < i_0 < i_1 < \dots$ is a strictly increasing sequence of indices, and let $(Q,T,Z) := \lim_{l \rightarrow \infty} (Q_{i_l},T_{i_l},Z_{i_l})$. It holds that $Q$ and $Z$ are unitary, and $T$ is upper triangular with a zero in position $(k,k)$. Moreover, 
    \begin{align*}
        Q T Z = \lim_{l \rightarrow \infty} Q_{i_l} \lim_{l \rightarrow \infty} T_{i_l} \lim_{l \rightarrow \infty}Z_{i_l} \stackrel{\text{(a)}}{=} \lim_{l \rightarrow \infty} Q_{i_l} T_{i_l} Z_{i_l} = \lim_{l \rightarrow \infty} V_{i_l} \stackrel{\text{(b)}}{=} U, 
    \end{align*}
where we used the product rule for finite limits in equality (a), and the fact that every subsequence of a convergent sequence has the same limit in equality (b). Hence, $U \in \text{Schur}^k_n$, and $\text{Schur}^k_n$ contains its accumulation points.
\end{proof}
To state Corollary \ref{cor:si}, let $\mathcal{Q}_n^k$ be the set of singular pencils having one right minimal index equal to $k$. (Note that pencils in $\mathcal{Q}_n^k$ may have more than one right minimal index.)
\begin{corollary}\label{cor:si}
    $\text{Schur}_n^{k+1} = \overline{\mathcal{S}_n^k} = \overline{\mathcal{Q}_n^k}$.
\end{corollary}
\begin{proof} It is plain by definition that $\mathcal{S}_n^k \subseteq \mathcal{Q}_n^k$. Let $U \in \text{Schur}_n^{k+1}$ and $T$ be a generalized Schur form of $U$ with $T_{k+1,k+1} = 0$. Theorem \ref{thm:arbitrarily_close} applied to $T$ shows that for every $\epsilon >0$ there is a neighbourhood of $U$ with radius $\epsilon$ having nonempty intersection with $\mathcal{S}_n^k$. Hence, $\text{Schur}_n^{k+1} \subseteq \overline{\mathcal{S}_n^k} \subseteq \overline{\mathcal{Q}_n^k}$. On the other hand, Lemma \ref{lemm:auxxk1} shows $\mathcal{S}_n^k \subseteq \mathcal{Q}_n^k  \subseteq \text{Schur}_n^{k+1}$ and Theorem \ref{thm:closure} shows $\text{Schur}_n^{k+1} = \overline{\text{Schur}_n^{k+1}}$. But, by definition, the closure of a set $X$ is the smallest closed set containing $X$, and hence, $\overline{\mathcal{S}_n^k}  \subseteq \overline{\mathcal{Q}_n^k} \subseteq \text{Schur}_n^{k+1}$.
\end{proof}

Having established that the minimization problem \eqref{eq:riproblem_k} is well posed, we now discuss how to solve it via Riemannian optimization. Computing a minimizer for \eqref{eq:riproblem_k} works analogously to finding a minimizer for \eqref{eq:riproblem}: we find a nearest upper triangular singular pencil with zero in a prescribed position, and use this knowledge to reformulate the problem in terms of a different objective function.

\begin{proposition} \label{prop:nearest_tri_k}
Let $A + \lambda B \in \C[\lambda]_1^{n \times n}$. An upper triangular singular pencil with zero in position $(k,k)$ nearest to $A+\lambda B$ is $\mathcal{P}_k(A) + \lambda \mathcal{P}_k(B)$ where
\[ \mathcal{P}_k(A)_{ij}=\begin{cases}A_{ij} \ &\mathrm{if} \ i<j \ \mathrm{or} \  i=j\neq k;\\
0 \ &\mathrm{otherwise}; \end{cases} \qquad \mathcal{P}_k(B)_{ij}=\begin{cases}B_{ij} \ &\mathrm{if} \ i<j \ \mathrm{or} \  i=j\neq k;\\
0 \ &\mathrm{otherwise}. \end{cases}\]
In particular, the squared distance to $A+\lambda B$ from $\mathcal{P}_k(A)+\lambda \mathcal{P}_k(B)$ is 
\begin{equation}\label{eq:objectivefunction_k} \mathcal{F}_k(A+\lambda B)= \sum_{i>j} (|A_{ij}|^2+|B_{ij}|^2) + (|A_{kk}|^2+|B_{kk}|^2 ).\end{equation}
\end{proposition}
\begin{proof}
Analogous to the proof of Proposition~\ref{prop:nearest_tri}.
\end{proof}

We consider the objective function 
\begin{equation}\label{eq:newobjfun_k}
     f_k(Q,Z)= \mathcal{F}_k(QAZ+\lambda QBZ) = [d(QAZ+\lambda QBZ,\mathcal{P}_k(QAZ)+\lambda \mathcal{P}_k(QBZ))]^2,  
\end{equation}
where $\mathcal{F}_k$ and $\mathcal{P}_k$ are as in Proposition~\ref{prop:nearest_tri_k} and $d$ is as in \eqref{eq:pencildist}. { Note that \eqref{eq:newobjfun_k} is equal to \eqref{eq:fkfirst}}. 

Computing \eqref{eq:riproblem_k} is equivalent to computing
\begin{equation}\label{eq:minnewobjfun_k}
    \min_{(Q,Z)\in U(n) \times U(n)} f_{k+1}(Q,Z) 
\end{equation} 
as well as an argument minimum $(Q_0,Z_0)$, which is stated more formally in Theorem \ref{thm:obj_function_k}.
\begin{theorem} \label{thm:obj_function_k}
Let $A+\lambda B \in \C[\lambda]_1^{n \times n}$. Let $f_{k+1}(Q,Z)$ be the function on $U(n) \times U(n)$ defined by \eqref{eq:newobjfun_k}. Then:
\begin{enumerate}
    \item The optimization problems \eqref{eq:riproblem_k} and \eqref{eq:minnewobjfun_k} have the same minimum value;
    \item The pair of unitary matrices $(Q_0,Z_0)$ is a global (resp. local) minimizer for \eqref{eq:minnewobjfun_k} if and only if the pencil $Q_0^* \mathcal{P}_{k+1}(Q_0 A Z_0) Z_0^* + \lambda Q_0^* \mathcal{P}_{k+1}(Q_0 B Z_0) Z_0^*$ is a global (resp. local) minimizer for \eqref{eq:riproblem_k}.
\end{enumerate}
\end{theorem}
\begin{proof}
The proof is analogous to that of Theorem~\ref{thm:obj_function}. Now, $\text{Schur}_n^{k+1}$ plays the role of $\mathcal{S}_n$; the set of upper triangular pencils of size $n\times n$ with a zero in position
$(k+1, k+1)$, $\text{Upp}_n^{k+1}$, plays the role of $\mathcal{T}_n$; and $\mathcal{F}_{k+1}$ plays the role of $\mathcal{F}$.
\end{proof}

Again, \eqref{eq:minnewobjfun_k} can be computed with Manopt \cite{manopt}; { note that, for all $k$, $f_{k+1} (Q,Z)$ are polynomial (hence everywhere smooth) functions in the entries of $Q$ and $Z$.} The Euclidean gradient and the Euclidean Hessian of $f_{k+1}$ are almost identical to those of $f$ (computed in Section~\ref{subs:complexgrad} and Section~\ref{subs:complexhess}, respectively); we only need to redefine $L$ such that now $L(C) + \lambda L(D) = (C-\mathcal{P}_{k+1}(C)) + \lambda ( D - \mathcal{P}_{k+1}(D))$, { which now is linear instead of piecewise linear}. The Riemannian gradient and the Riemannian Hessian are obtained from their Euclidean counterparts as explained in Section \ref{subsec:riemanprojec}.

% \begin{vn}
%     I would personally remove both blue bits above. No need to write once more that (3.3) was not smooth, we said it already. The linear vs piecewise linear also is, in my opinion, more of a distraction than a useful comment.
% \end{vn}

% \begin{fd}
%     Vanni, please proceed as you prefer. No problem on my side. In my view, the second blue sentence at least is useful for the reader, taking into account that a similar comment was made in Section 3.
% \end{fd}

% \begin{vn}
%     I am OK about the content of Remark 5.9 (I only edited slightly), but I would prefer to say it in the response than in the paper, simply because the observation was a bit casual: when I looked at MATLAB, this seemed to be the case. We never did e.g. statistical experiments to confirm that this behaviour is typical. Even if it is probably correct, I am not so fully confident to state it in the article (but I think it would be OK to mention it in the response).
% \end{vn}

% \begin{fd}
%     I removed the previous Remark 5.9. I agree that the right way to proceed is not including statements about which we are not fully sure. They do not help to make the paper better or more solid. Also, I think that there is no need to mention anything about this remark in the response letter.
% \end{fd}

{\section{A Riemannian algorithm for the nearest real singular pencil} \label{sec:real}
In applications, sometimes it is desirable to add the real structure to the problem, that is, $A$ and $B$ are real matrices and one wants to find the singular pencil $S + \lambda T$ nearest to $A + \lambda B$ where $S,T$ are also real. Temptingly, one could imagine to simply replace the special unitary group with the special orthogonal group, without otherwise changing the analysis, or the algorithms in Sections \ref{sec:complex}, \ref{sec:smoothalgorithm} and \ref{sec:third}, except for the trivial modifications that take into account that all the matrices are now real and not complex. However, an apparent difficulty is that Lemma \ref{lemma:schur} does not hold for a general real pencil unless one relaxes the upper-triangularity constraint to only being block  upper triangular, having possibly $2 \times 2$ real blocks on the diagonal. Crucially, though, the proof of Theorem \ref{thm:obj_function} only requires the fact that every \emph{singular} pencil has a (properly triangular) Schur form. The following Lemma \ref{lem:froilannew} demonstrates that this is indeed true also over the real field.
\begin{lemma}\label{lem:froilannew}
    Let $S + \lambda T \in \R [\lambda]_1^{n \times n}$ be a real singular pencil, that is, $\det (S + \lambda T) \equiv 0$. Then, there exist $Q,Z \in SO(n)$ such that $QSZ$ and $QTZ$ are both (real and) upper triangular.
\end{lemma}
\begin{proof}
    Let $X,Y \in GL(n,\R)$ be such that  $K  = X^{-1}(S + \lambda T)Y^{-1} \in \R[\lambda]_1^{n \times n}$ is a real Kronecker canonical form of $S + \lambda T$, which is unique up to permutation of its diagonal blocks: see \cite[Theorem 3.2]{real_KCF}. Without loss of generality, we may assume that $K=K_1 \oplus K_2$ where $K_1$ contains all, and only, the blocks associated with the right minimal indices whereas $\ker K_2 = \{ 0 \}$. Then, by construction, $K$ is upper triangular. Consider now the QR decomposition $X=Q_0^T R_X$ and the RQ decomposition $Y=R_Y Z_0^T$. Then, $Q_0 (S + \lambda T) Z_0 = R_X K R_Y$ is a triangular pencil, and hence $Q_0 S Z_0$ and $Q_0 T Z_0$ are both upper triangular with $Q_0,Z_0 \in O(n)$. Similarly to Lemma \ref{lem:SUnotU}, it now suffices to define $Q=(I_{n-1} \oplus \det Q_0) Q_0 \in SO(n)$ and $Z=Z_0 (I_{n-1} \oplus \det Z_0)  \in SO(n)$ to conclude the proof.
\end{proof}

An immediate consequence of Lemma \ref{lem:froilannew} is that the, apparently problematic, idea of simply running our algorithms in real arithmetic and over the (connected!) manifold $SO(n) \times SO(n)$ actually works. We omit a detailed analysis as it is completely equivalent, mutatis mutandis, to what we have carried out for the complex case in Sections \ref{sec:complex}, \ref{sec:smoothalgorithm} and \ref{sec:third}.
}

\section{Implementation}\label{sec:imp}

We minimize the { appropriate} objective function numerically by using the trust-region solver \cite{ABG07} of Manopt \cite{manopt} version 7.1. It uses a descent method that finds the desired update direction in an inner iteration via a Newton-type quadratic subproblem. This subproblem is solved via a truncated conjugate-gradient method outlined in \cite[Algorithm~11]{ABG07}. An outline for an algorithm implementing the trust-regions solver is given in \cite[Algorithm~10]{ABG07}.

Manopt uses the norm of the gradient as a measure of convergence, and stops the iteration once the gradient norm falls below a prescribed (absolute) threshold. In order to make this stopping criterion invariant to scaling, in our implementation we normalize the input to have a fixed norm. However, the implementation of the trust-regions method in Manopt is highly sophisticated and not entirely independent of the scaling (even after the stopping criterion threshold is scaled accordingly), and hence the value of the input norm can have an effect on the performance. Heuristically, a choice that seems to work well in practice consists of normalizing the pencil to have norm $100$ and setting the parameters within Manopt so that the iteration stops once the gradient norm evaluates to less than $10^{-10}$; these are therefore the default values in our implementation.  After the scaled problem is solved, we multiply the obtained singular pencil by the inverse of the introduced scale factor to get the final solution. Manopt also requires as input a starting point on the manifold $U(n) \times U(n)$. We discuss in Section \ref{sec:numerical} several strategies for choosing the starting point.

{ MATLAB codes that implement the non-smooth algorithm in Section \ref{sec:complex}, the smoothed algorithm in Section \ref{sec:smoothalgorithm} and the algorithm in Section \ref{sec:fixminindex} for the problem with a prescribed minimal index are} publicly available in the {GitHub} repository at \href{https://github.com/NymanLauri/nearest-singular}{https://github.com/NymanLauri/nearest-singular}\href{https://github.com/}. { These codes require} the Manopt package, which can be downloaded from \href{https://www.manopt.org/}{manopt.org}.

{
\begin{remark}
    Note that, for the cost function \eqref{eq:fqz}, the Euclidean gradient, the Euclidean Hessian, and their projections onto the tangent space $T_Q U(n) \times T_Z U(n)$ (presented in Section \ref{sec:complex}) are all real for a real pencil $A + \lambda B$, real arguments $(Q,Z)$, and real directions. Moreover, Manopt performs retraction for real inputs using a QR-decomposition which yields a real orthogonal matrix; a detailed description on retractions can be found in \cite[Sections 7.3 and 7.4]{rmanifolds}. These facts combined guarantee that the algorithm in Section \ref{sec:complex} for the complex version of the problem gives a real output whenever the  input pencil and the starting point $(Q_0, Z_0)$ are real. Taking into account the discussion in Section \ref{sec:real}, this means that the real algorithm is incorporated in a natural way in the complex one just by using real inputs and real starting points.
\end{remark}
}

\section{Numerical experiments}
\label{sec:numerical}
We present several numerical experiments within this section. Unless stated otherwise, we consider the standard algorithm presented in Section~\ref{sec:complex} implemented as in Section \ref{sec:imp}. The algorithm that allows us to prescribe the minimal index of the sought singular pencil, presented in Section \ref{sec:fixminindex}, is used solely in Subsection~\ref{sec:minindex_numerical}{, while the smoothed algorithm outlined in Section \ref{sec:smoothalgorithm} is used solely in Example \ref{ex71} and Subsection \ref{sec:smooth_numerical}.}

The structure of this section is as follows. Subsection~\ref{sec:examples} shows how our method performs with a selection of example pencils that were proposed, and tackled by alternative methods, in \cite{bhm,bora,glm}. In Subsection~\ref{sec:running_time}, we characterize the running time of our code, in particular how the running time grows in the size of the pencil. We examine in Subsection~\ref{sec:initial} the effect of the starting point on the output.  In Subsection~\ref{sec:multiple_starting}, we measure how much the solution can be improved by running the algorithm multiple times from different starting points, and characterize the distribution of local minima for randomly generated  input pencils. Subsection~\ref{sec:comparison} outlines a statistical experiment that compares our method with the approach of \cite{bora}, which is arguably the best existing algorithm prior to this paper. We also compare in part with \cite{glm}, although in the latter case the lack of an efficient implementation of 
 the algorithm in \cite{glm} limits heavily the size of the input for a statistical experiment. {In Subsection~\ref{sec:minindex_numerical}, we demonstrate the possibility of prescribing the minimal index, and argue that this can help reveal meaningful information about the problem. Finally, in Subsection \ref{sec:smooth_numerical} we compare the performance of the smoothed (Section \ref{sec:smoothalgorithm}) and non-smoothed  (Section \ref{sec:complex}) algorithms.}

\subsection{Individual examples} \label{sec:examples}

\begin{example}\label{ex71}
    A model of a two-dimensional, three-link mobile manipulator is derived in \cite{real_application}. This model generates an $8\times 8$ matrix pencil $A + \lambda B$ such that
    \begin{equation}\label{eq:difficultpencil}
         A=\left[\begin{array}{ccc}
    0 & I_3 & 0 \\
    -K_0 & -D_0 & F_0^{\mathrm{T}} \\
    F_0 & 0 & 0
    \end{array}\right], \quad B=\left[\begin{array}{ccc}
    I_3 & 0 & 0 \\
    0 & M_0 & 0 \\
    0 & 0 & 0
    \end{array}\right],
    \end{equation}
   
    where
    \begin{align*}
    M_0 & =\left[\begin{array}{rrr}
    18.7532 & -7.94493 & 7.94494 \\
    -7.94493 & 31.8182 & -26.8182 \\
    7.94494 & -26.8182 & 26.8182
    \end{array}\right], D_0=\left[\begin{array}{rrr}
    -1.52143 & -1.55168 & 1.55168 \\
    3.22064 & 3.28467 & -3.28467 \\
    -3.22064 & -3.28467 & 3.28467
    \end{array}\right], \\
    K_0 & =\left[\begin{array}{rrr}
    67.4894 & 69.2393 & -69.2393 \\
    69.8124 & 1.68624 & -1.68617 \\
    -69.8123 & -1.68617 & -68.2707
    \end{array}\right], F_0=\left[\begin{array}{lll}
    1 & 0 & 0 \\
    0 & 0 & 1
    \end{array}\right] .
    \end{align*}
For the above pencil, the { Riemannian algorithm takes an unusually long time to converge both in the non-smooth and smoothed variants} (see Subsection~\ref{sec:minindex_numerical} for a potential cause of this behaviour). Eventually, it finds a singular pencil at a distance of 0.01117 (this pencil is stored at \href{https://github.com/NymanLauri/nearest-singular}{https://github.com/NymanLauri/nearest-singular} as \texttt{example6.1.mat}). This agrees with the result given by the ODE approach \cite[Example 2.20]{glm} and the approach in \cite{bora} (at least taking into account that both these two sources report fewer significant digits than us), as well as the lowest upper bound in \cite[Example 14]{bhm}. 
\end{example}

\begin{example}
Consider the $3 \times 3$ matrix pencil
$$
A + \lambda B=\left[\begin{array}{lll}
1 & 0 & 0 \\
0 & \epsilon & 0 \\
0 & 0 & 1
\end{array}\right]-\lambda\left[\begin{array}{lll}
0 & 1 & 0 \\
0 & 0 & 1 \\
0 & 0 & 0
\end{array}\right].
$$    
For this example, the distance to singularity can be computed to be $\epsilon$ \cite[Example 4]{bhm}. We fix $\epsilon = 10^{-8}$ and run our algorithm starting from different initial points. If the Cartesian product of two identity matrices is used as the initial point, the algorithm converges in just one iteration to the global minimum $\epsilon$. Starting from a random point, the algorithm sometimes computes $\epsilon$ which is the global minimum, and sometimes computes 1 which is a local minimum. After running the algorithm $10^{4}$ times starting from different random points, the global minimum $\epsilon$ was computed $38.81 \%$ of the time, whereas the local minimum 1 was computed $61.19 \%$ of the time. %

\end{example}

\begin{example}
Let us consider a pencil of form $B_n - \lambda B_n$, where $B_n \in \mathbb{R}^{n \times n}$ is defined by
$$
B_n=\left[\begin{array}{ccccc}
1 & -1 & -1 & \ldots & -1 \\
& 1 & -1 & \ldots & -1 \\
& & \ddots & \vdots & \vdots \\
& & & 1 & -1 \\
& & & & 1
\end{array}\right].
$$
For $n=20$, a lower bound for the distance to singularity is $4 \times 10^{-6}$ \cite[Example 5 in Table 2]{bhm}. Starting from a random point, our algorithm consistently obtains this result and hence converges to a global minimum.
\end{example}

\begin{example}
Let us consider the pencil
\begin{align} \label{eq:bhm_ex_6}
A+\lambda B=\left[\begin{array}{cc}
1 & \frac{1}{\epsilon} \\
0 & 1
\end{array}\right]-\lambda\left[\begin{array}{cc}
0 & \frac{1}{\epsilon} \\
0 & 1
\end{array}\right]
,
\end{align}
where $\epsilon > 0$. For a wide range of values for $\epsilon$ ($\epsilon=10^{-2}, 10^{-4}, 10^{-6}, 10^{-8}, 10^{-10}, 10^{-12},$ $10^{-14}$), the best result given by our algorithm is consistently equal to $\epsilon$, which is a lower bound according to \cite[Example 6 in Table 2]{bhm} and hence a global minimum.  The authors of \cite{bora} also report finding the global minimum $\epsilon$ \cite[Example 9.5]{bora} for $\epsilon=10^{-4}$.

\end{example}

\begin{example}
    In \cite[Example 9.6 and Example 9.7]{bora}  as well as in \cite[Example 2.15 and Example 2.16]{glm}, the pencils
    \[\begin{bmatrix}0&0.04&0.89\\
        0.15&-0.02&\lambda\\
        0.92&\lambda+0.11&0.066
    \end{bmatrix} \ \mathrm{and} \ \begin{bmatrix}
        -1.79&0.1&-0.6\\
        0.84&-0.54&\lambda+0.49\\
        -0.89&\lambda+0.3&0.74
    \end{bmatrix} \]
    are considered. Our algorithm computes the same potential global minima proposed in \cite{bora} with, respectively, a distance $\approx 0.1155462894$ with the nearest singular pencil having a right minimal index $1$ for the first input, and a distance $\approx 0.9435641675$ with the nearest singular pencil having a right minimal index $2$ for the second input. These computed distances are lower than the values $0.1193$ and $1.065$, respectively, reported in \cite[Example 2.15 and Example 2.16]{glm}. Note that we can find out the right minimal index of each computed singular pencil as a consequence of Theorem \ref{thm:arbitrarily_close}.
\end{example}

\subsection{Running time} \label{sec:running_time}
\begin{figure}[h]
\begin{subfigure}{.5\textwidth}
  \centering
  \includegraphics[width=\linewidth]{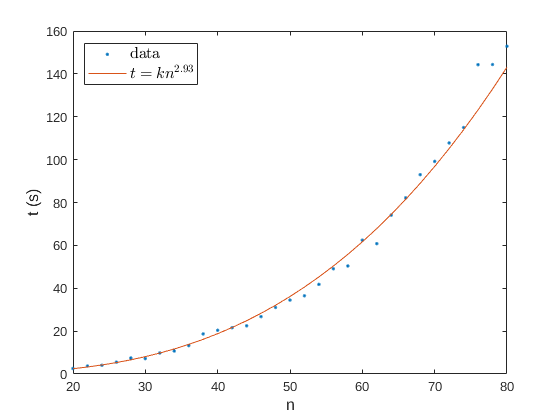}
  \caption{Interval $n \in [20,80]$ }  
  \label{fig:running_time_small}
\end{subfigure}%
\begin{subfigure}{.5\textwidth}
  \centering
  \includegraphics[width=\linewidth]{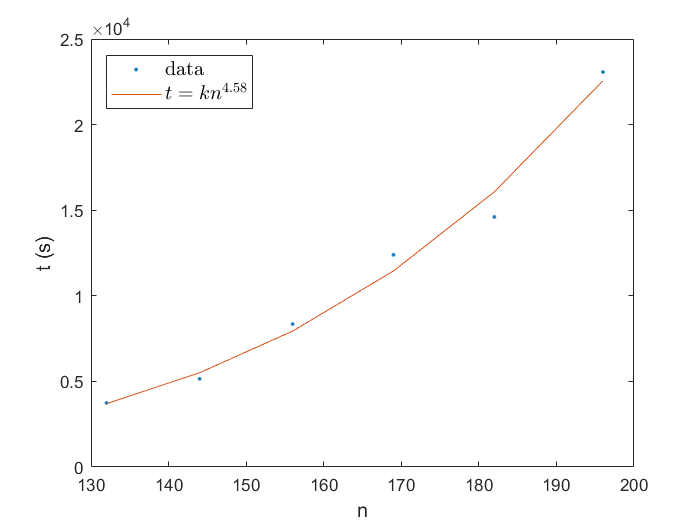}
  \caption{Interval $n \in [130,200]$ }
\label{fig:running_time}
\end{subfigure}
\caption{Average running time versus size. For the interval $n \in [20,80]$, the least squares fit yields approximately $t = k\, n^{2.93}$, where $k \approx 3.8310 \times 10^{-4}$, whereas for the interval $n \in [130,200]$, the least squares fit yields approximately $t = k\, n^{4.58}$, where $k \approx 7.3423 \times 10^{-7}$.}
\label{fig:running_times}
\end{figure}

For each $n$, we generate random $n \times n$ pencils $A + \lambda B$ s.t. the real and imaginary part of each element in $A$ and $B$ is drawn from the normal distribution $\mathcal{N}(0,1)$. We measure how long it takes for the algorithm to fully converge. This was repeated 50 times for each $n \in [20,80]$ after which the average running times were computed. The computation was run using MATLAB R2023a on an Intel Core i5-12600K, a processor commonly found in ordinary computers. Finally, a least squares fit was made for the loglog-plot in order to estimate the order of the growth rate in $n$. Figure~\ref{fig:running_time_small} shows that the running time in the interval $n \in [20,80]$ approximately follows the line $t = k\, n^{2.93}$, where $k \approx 3.8310 \times 10^{-4}$. Hence, increasing $n$ increases the running time relatively mildly in this interval. Indeed, the average running time for $n=80$ was under three minutes! We stress that this is an enormous breakthrough with respect to other existing methods, that on the same machine seemed to typically take much longer than a few minutes to converge already for smaller inputs.

We have run a similar experiment for larger sizes in the interval $n \in [130,200]$ with 24 runs for each $n$. For this more demanding experiment, we used MATLAB's internal parallelization with 24 processes on a 2x12 core Xeon E5 2690 v3 2.60GHz. Figure~\ref{fig:running_time} shows that in this case the increase is much steeper, and suggests $O(n^{4.58})$ as the asymptotic behaviour.

\subsection{Initial guess} \label{sec:initial}
As the problem is non-convex, the choice of the initial guess can alter the local minimum to which the method converges. Hence, an educated guess for an initial point can potentially improve the performance of the algorithm. To see whether we observe this in practice, we compare the performance of different initial guesses against each other. In particular, we consider the following three types of initial guesses:
\begin{enumerate}[{(a)}]
    \item Random point. 
    \item The point that gives minimal value for the objective function among 10 randomly generated points.  
    \item An initial guess that brings the pencil into Schur form. More precisely, we consider different cyclic permutations of the eigenvalues in the Schur form, and choose the one that minimizes the objective function. 
\end{enumerate}

We randomly generated $10^4$ complex $10 \times 10$ pencils $A + \lambda B$, by drawing the real and imaginary part of each element in $A$ and $B$ from the normal distribution $\mathcal{N}(0,1)$. For each of them we generated one initial guess of each of the three types described above, and ran the algorithm once for each type of initial points. The random points for (a) and (b) were constructed by drawing $Q$ and $Z$ at random from the Haar measure on $U(n)$; the same applies throughout this section whenever we refer to ``random initial points" or similar phrases. For the initial guess (c), we first compute one such point $(Q,Z)$ via the \texttt{qz}-function in MATLAB. Then, we go through all cyclic permutations of eigenvalues of the form $(1\ 2\ 3\ \dots \ k) (k+1) (k+2)...(n)$, $k=2,3,\dots,n$, where $n$ is the size of the pencil. To be more precise, we consider Schur forms which result from successively applying transpositions $(1\ 2), \dots , (k-1 \ k)$ on the eigenvalues (see Theorem~\ref{th.11} for the details on swapping two eigenvalues). This is repeated for $k=2,3,\dots,n$.

\begin{table}[h]
\begin{center}
\begin{tabular}{cccc} \toprule 
    {Initial guess} & {Frequency of best output} & {Median distance} & {Average distance}  \\ \midrule
    {(a)}  & 33.74 \%  & 1.8059 & 1.8767  \\
    {(b)}  & 33.47 \%  & 1.8059 & 1.8763  \\ 
    {(c)}  & 32.79 \%  & 1.7914 & 1.8332  \\ \bottomrule
\end{tabular}
\caption {Comparison of different initial points} \label{tab:initial_guess}
\end{center}
\end{table}

In Table~\ref{tab:initial_guess}, it can be seen that all the three initial points show almost identical behaviour in terms of the median and average distances and how often that particular choice yielded the best solution. For the average of $10^4$ Bernoulli distributed random variables with mean $\frac{1}{3}$, the standard deviation is $\frac{\sqrt{2}}{3} \cdot 10^{-2} \approx 0.5 \%$. Hence, the difference between each initial guess and the mean $\frac{1}{3}$ is approximately one sigma. Consequently, the choice between these three types of initial points does not seem to significantly affect the performance of the algorithm.

\subsection{Multiple starting points and distribution of local minima} \label{sec:multiple_starting}

We saw in Figure~\ref{fig:running_times} that the running time of our algorithm is reasonably fast, e.g., our code usually converges within a few minutes for randomly generated input pencils of size $n \lesssim 100$. Since this is a tremendous speedup with respect to competitor methods, it becomes much less time-consuming to compute several local minima starting from different initial points, and eventually select the best answer. (Moreover, note that this task can be trivially parallelized.) In the following experiment, we try to characterize how much the solution improves when the algorithm is run multiple times. 

We randomly generate a $5\times5$ pencil s.t. the real and imaginary parts are drawn from $\mathcal{N}(0,1)$. For different values of $n$, we run the algorithm $n$ times from different random starting points. We compare the best output with the result of the first run, and compute the ratio of these distances, which we denote by $\delta_r$. We repeat this $10^3$ times, and compute the average value for the relative distance $\delta_r$ which we denote by $\bar \delta_r$. Finally, we repeat this experiment for randomly generated $10\times10$ and $20\times20$ pencils.

\begin{figure}[h]
  \centering
  \includegraphics[width=0.6\linewidth]{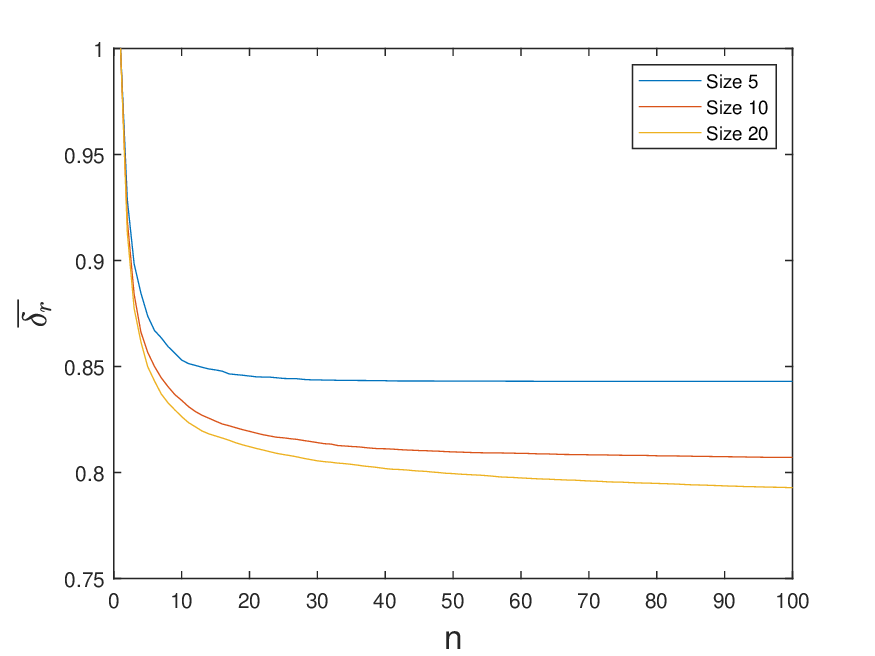}
\caption{For $n \in [1,100]$, the distance to the best output amongst $n$ runs relative to the distance to the initial output.}
\label{fig:improvement}
\end{figure}

Figure~\ref{fig:improvement} shows that the ratio between the distance to the best output and the distance to the initial output first decreases rapidly and eventually settles at a constant value as $n$ grows large; this is of course expected in the sense that for very large values of $n$ one would eventually hit a starting point from which the method actually converges to the global minimum. While this general trend is the same for all sizes of the pencil, larger sizes exhibit a more significant improvement in the distance. Table~\ref{tab:improvement} shows the exact values for this improvement for a selection of values for $n$. For example, for $20 \times 20$ pencils, the experiment suggests that one can expect to improve the distance by roughly $15 \%$ by simply running the algorithm $5$ times from $5$ different starting points. 

\begin{figure}[h]
  \subcaptionbox{Results for one randomly generated $10\times10$ pencil. The computed $10^2$ minima correspond to  46 distinct pencils.}%
  {\includegraphics[width=0.48\linewidth]{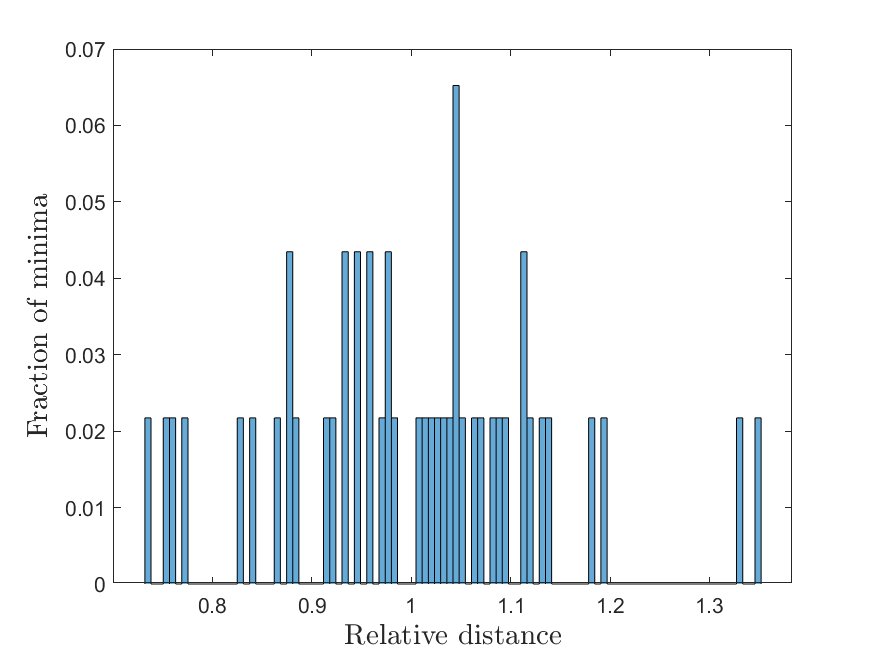}}
  \hspace{\fill}
  \subcaptionbox{Aggregate over $10^3$ randomly generated $10\times10$ pencils. On average, 33.8 of the computed $10^2$ singular pencils are distinct.}%
  {\includegraphics[width=0.48\linewidth]{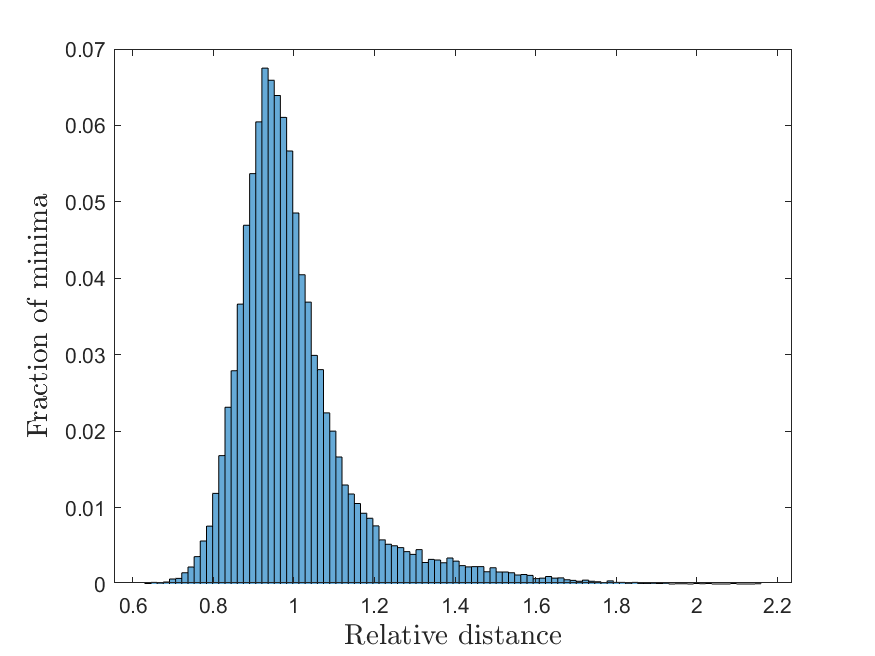}}
  \hspace{\fill}
  \subcaptionbox{Results for one randomly generated $20\times 20$ pencil. The computed $10^2$ minima correspond to 63 distinct pencils.}%
  {\includegraphics[width=0.48\linewidth]{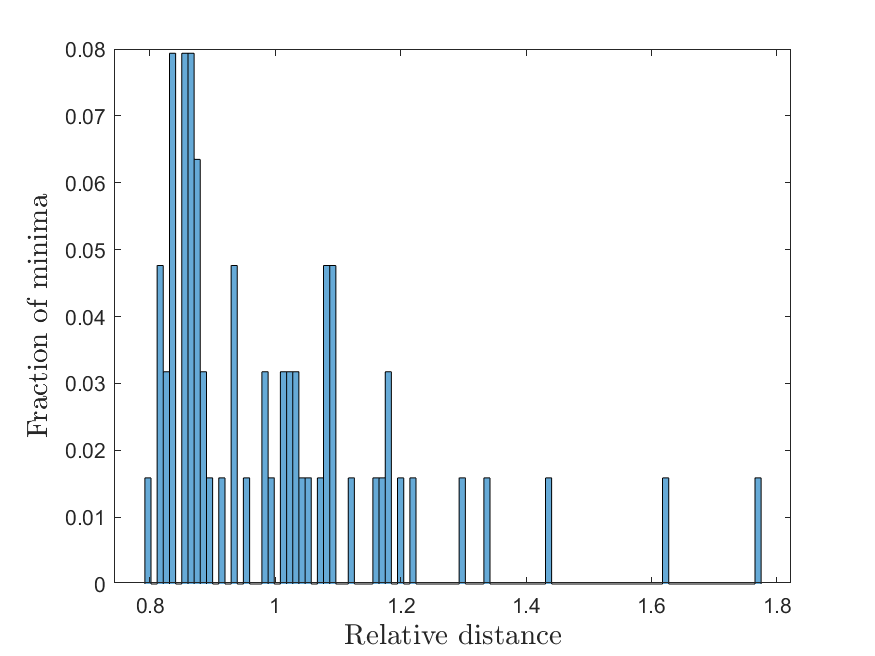}}
    \hspace{\fill}
  \subcaptionbox{Aggregate over $10^3$ random $20\times20$ pencils. On average, 70.3 of the computed $10^2$ singular pencils are distinct.}%
  {\includegraphics[width=0.48\linewidth]{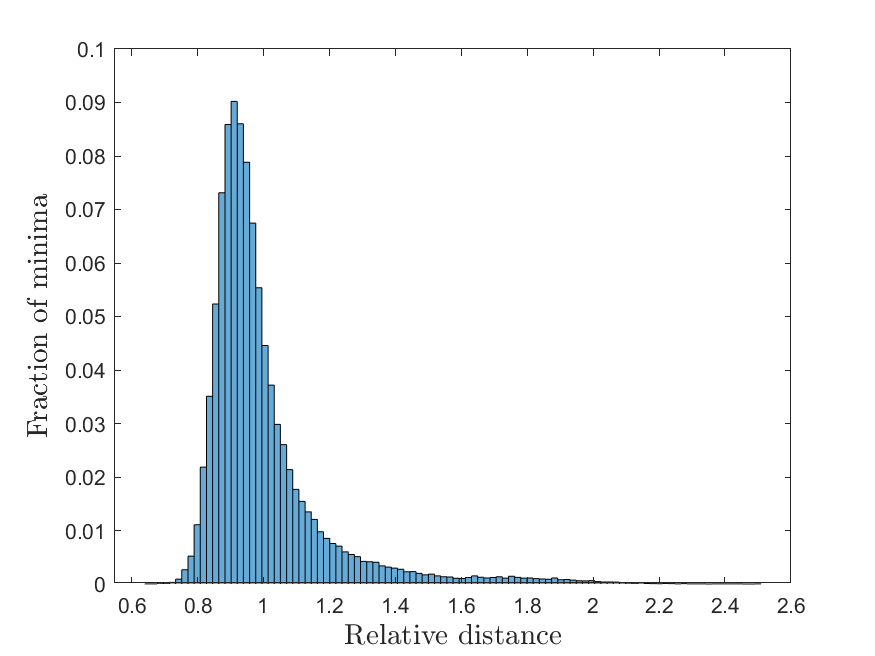}}

  \caption{The distribution of local minima for randomly generated $10\times 10$ and $20\times20$ pencils, where the distance is given relative to the average of the local minima corresponding to the pencil.
For each pencil, $10^2$ minima were computed, each starting from a different random point.}\label{fig:distribution}
\end{figure}

This type of experiment can effectively provide us with some insight on the problem, and in particular on the distribution of the local minima. To this goal, we conduct a numerical experiment with the aim of characterizing the distribution of local minima for randomly generated $10\times10$ and $20\times20$ pencils. To do this, we first randomly generate a $10\times10$ pencil s.t. the real and imaginary parts are drawn from $\mathcal{N}(0,1)$. We then run the algorithm $10^2$ times from different random starting points, and record the values of the computed minima. It is possible that the same singular pencil was computed multiples times, and the probability of computing a certain output may depend not only on the problem but also on the particular method that is employed to solve it. Thus, we remove duplicates from the data by using the tolerance $10^{-8}$ (i.e. two values are considered the same if they differ from each other in the relative sense by less than the tolerance), and arrange the values of the computed local minima in a histogram such that the distances are given relative to the average. We repeat this $10^3$ times for different random inputs, and aggregate the relative distances in a single histogram. The same process is repeated for random $20\times20$ pencils. Figure~\ref{fig:distribution} shows that for both sizes, the computed singular pencils corresponding to local minima are more commonly located close to the average distance, while no computed singular pencil was at a distance of more than three times the average, or less than half the average. Thus, the experiments suggest that, although one cannot guarantee that a global minimum was computed, it is likely that the computed (local) minimum provides a reasonable estimate of the distance to singularity.

\begin{table}[h]
\begin{center}
\begin{tabular}{cccc} \toprule 
    {} & \multicolumn{3}{c}{Improvement}  \\ 
    {$n$} & Size 5 & Size 10 & Size 20  \\ \midrule
    {2}  & 7.14 \% & 8.12 \% & 8.71 \%  \\
    {3}  & 10.16 \% & 11.59 \% & 12.25 \%  \\ 
    {4}  & 11.54 \% & 13.38 \% & 13.83 \%  \\ 
    {5}  & 12.62 \% & 14.33 \% & 15.01 \%  \\ 
    {100}  & 15.70 \% & 19.29 \% & 20.72 \%  \\ \bottomrule
\end{tabular}
\caption {Change in the computed distance as the algorithm is run multiple times choosing the best output. Improvements are measured relatively to the initial output.} \label{tab:improvement}
\end{center}
\end{table}

\subsection{Comparison with other existing algorithms}\label{sec:comparison}
\begin{figure}[h]
  \subcaptionbox{Mean distance}%
  {\includegraphics[width=0.49\linewidth]{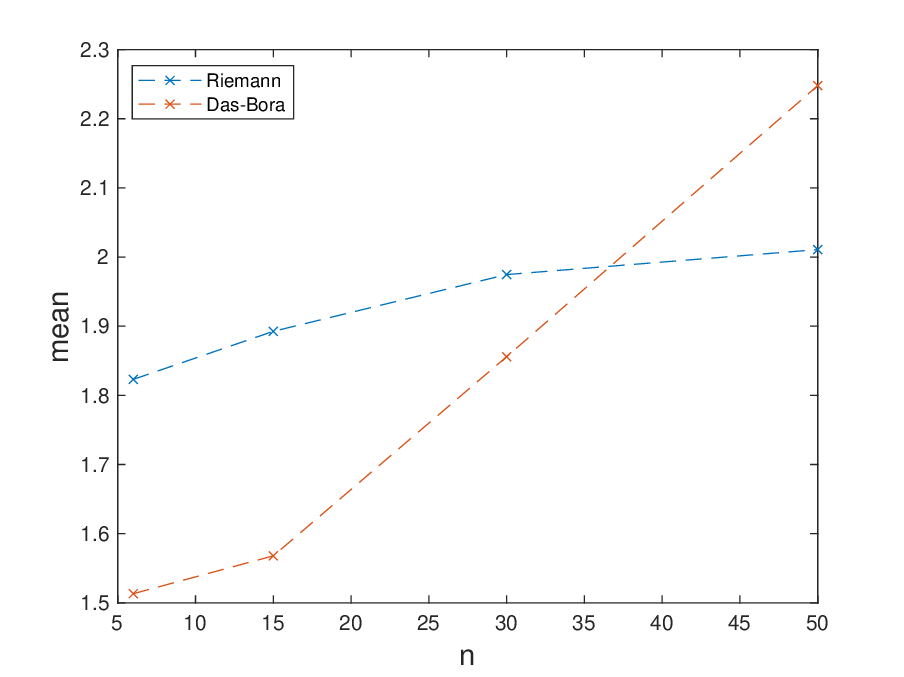}}
  \hspace{\fill}
  \subcaptionbox{Median distance}%
  {\includegraphics[width=0.49\linewidth]{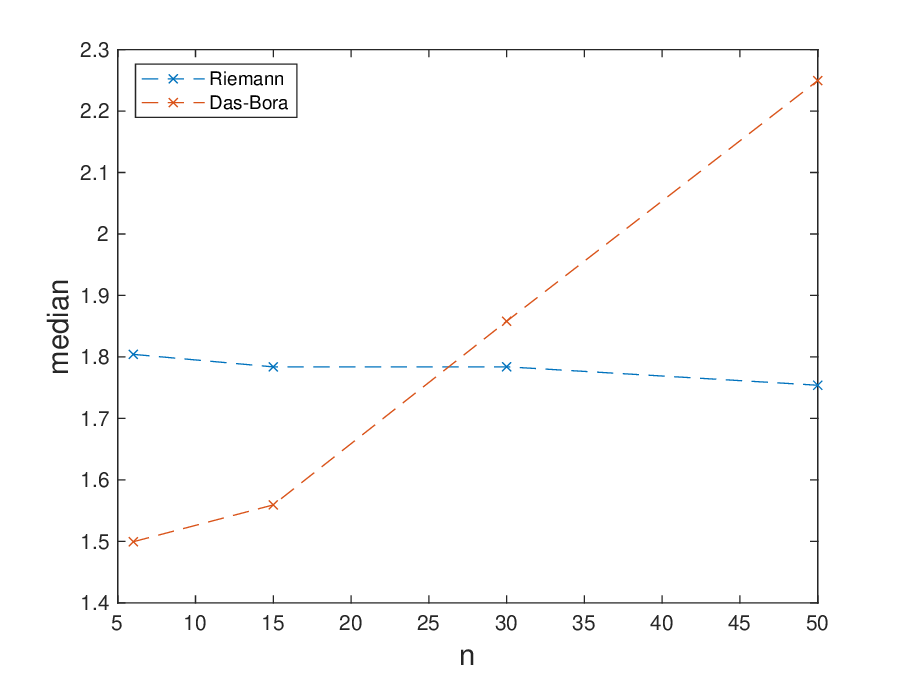}}
  \hspace{\fill}
  \subcaptionbox{Performance profile}%
  {\includegraphics[width=0.49\linewidth]{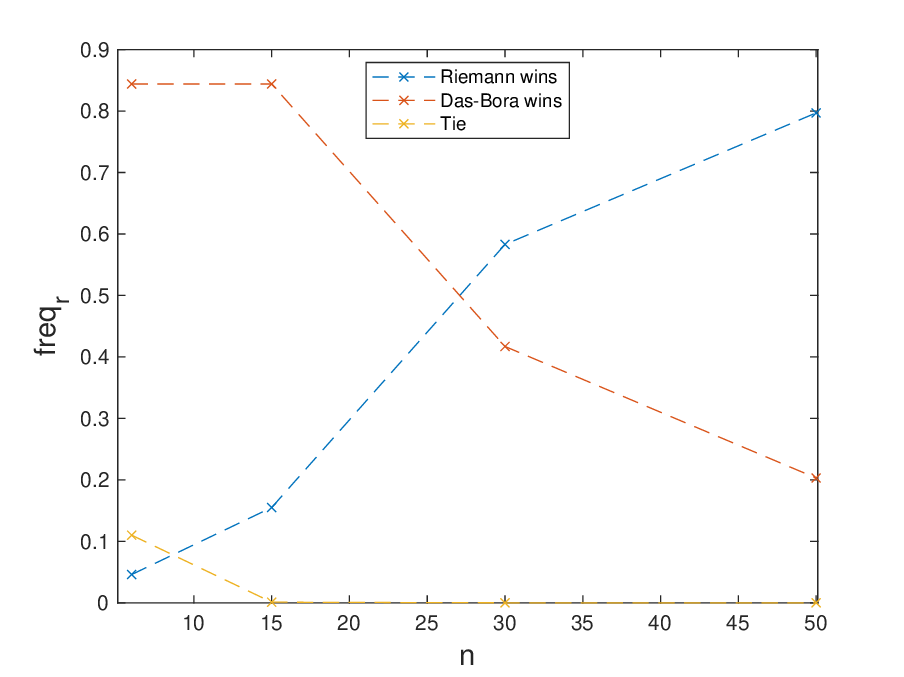}}
    \hspace{\fill}
  \subcaptionbox{Running time}%
  {\includegraphics[width=0.49\linewidth]{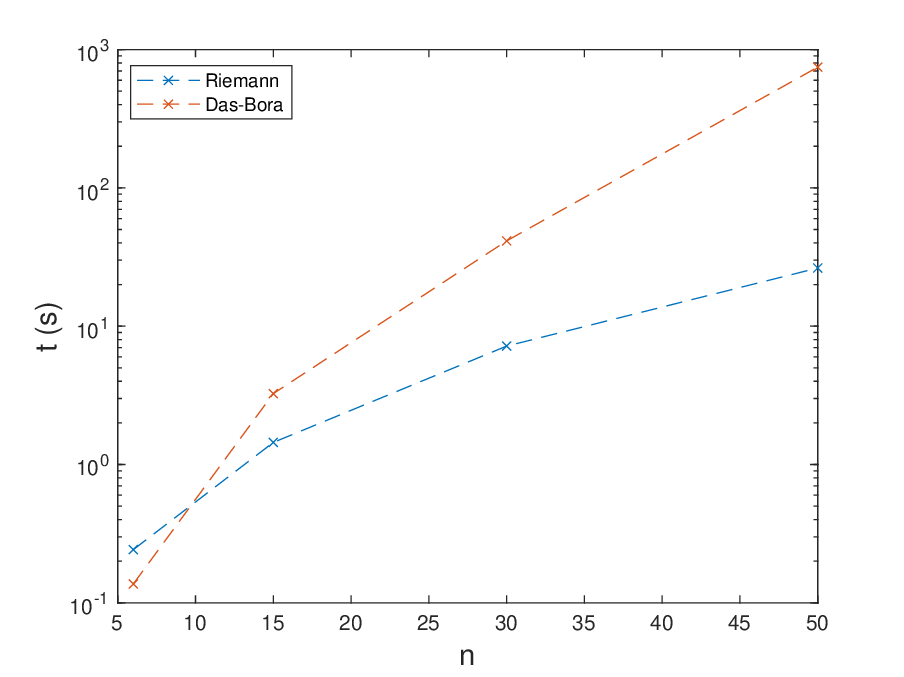}}

  \caption{Comparison between the method of this paper (Riemann) and the Das-Bora algorithm \cite{bora} for $n \in \{6,15,30,50\}$. The performance profile reports the relative frequency of which method yielded a better solution (or of ties), while running times were measured using MATLAB R2023a on an Intel Core i5-12600K.}\label{fig:bora_comparison}
\end{figure}

We benchmark the me\-thod presented in this paper against what we consider the two currently best methods for computing the nearest singular pencil, i.e., the very recent algorithm\footnote{For the algorithm in \cite{bora}, we have used the MATLAB codes kindly provided by the authors of \cite{bora}. In particular, we have used the algorithm in \cite{bora} based on BFGS, which is the fastest one in our experience.} of \cite{bora} and the ODE-approach\footnote{We implemented the algorithm in \cite{glm} based on the pseudo-code \cite[Algorithm 1]{glm}. For the inner iteration, we used the MATLAB codes kindly provided by the authors of \cite{glm}.} presented in \cite{glm}. We do this by comparing the output of the three
methods for complex randomly generated pencils of different sizes. Having fixed a size $n \in \{6,15,30,50\}$, we generated random $n \times n$ pencils $A + \lambda B$ s.t. the real and imaginary part of each element in $A$ and $B$ is drawn from the normal distribution $\mathcal{N}(0,1)$. After $10^3$ runs, we compared our algorithm (implemented by running it only once with  $(I_n,I_n)$ as the starting point, where $I_n$ denotes the $n\times n$ identity matrix) with \cite{bora} in terms of quality of the output (we report both the mean and the median, as well as the relative frequency of each method providing the best solution), and running time. The results of this comparison are reported in Figure \ref{fig:bora_comparison}. The quality of the output of our algorithm was typically worse than \cite{bora} for very small inputs $n=6,15$, but slighlty better for $n=30$ and already much better for $n=50$. In terms of running time, we outperformed \cite{bora} already for $n=15$; for $n=50$ the difference was already striking, with a ratio of average running times $\approx 29$ in favour of our method.

We next compare with \cite{glm} for $n=6$. The results of the comparison are seen in Table~\ref{tab:ode_comparison}. The approach detailed in this paper finds a better solution approximately 60 \% of the time, and computes an approximately 11 \% smaller distance on average. The amount of times the two algorithms converged to the same value is not recorded in the table as this happened a negligible amount of time. We note that \cite{bora}, while not included in Table~\ref{tab:ode_comparison}, often does better than both \cite{glm} and our method for the very small size $n=6$, as can be seen from the leftmost point of the graphs in Figure \ref{fig:bora_comparison}.

\begin{table}[h]
\begin{center}
\begin{tabular}{cccc} \toprule
    {Method} & {Frequency of best output} & {Median distance} & {Average distance}  \\ \midrule
    {ODE \cite{glm}}  & 37.3 \% & 1.8925 & 2.0601  \\
    {This paper}  & 62.7 \% & 1.8042  & 1.8231 \\ \bottomrule
\end{tabular}
\caption {Comparison with the ODE approach for $6 \times 6$ pencils.} \label{tab:ode_comparison}
\end{center}
\end{table}

In summary, the experiments suggest that \cite{bora} is usually the best method for extremely small inputs, but as the size increases our novel approach soon dominates. This is not surprising since our Riemannian algorithm, unlike \cite{bora} that can deal with any polynomial matrix, is specialized to pencils. We also note that, in the comparison, we used the most naive (one-run) implementation of our algorithm: as the running time of the Riemannian method quickly becomes much faster than \cite{bora} and as discussed in Subsection \ref{sec:multiple_starting}, we might also have run it multiple times starting from different random initial points while keeping the computational time below that of \cite{bora}, this would have decreased further the threshold where the approach of this paper starts doing better than \cite{bora} in terms of the quality of the computed minimizers.  

It would be interesting to extend the comparison of the three methods to pencils of even larger size. However, this is impractical because the ODE-algorithm implemented based on the pseudo-code \cite[Algorithm 1]{glm} is too slow to perform an extensive statistical analysis within a reasonable amount of time for $n > 6$. The same comment holds for the algorithm of \cite{bora} when $n>50$. Extrapolating from the trend of Figure \ref{fig:bora_comparison}, we conjecture that our method would asymptotically dominate over \cite{bora} not only for speed but also for output quality. Regarding \cite{glm}, at least until a more efficient implementation is developed, its applicability is in practice limited to very small sizes only.

\subsection{Nearest pencil with prescribed minimal index} \label{sec:minindex_numerical}

We revisit Example \ref{ex71} by computing the distances to the nearest pencil with a prescribed right minimal index $k$, $0 \leq k \leq 7$. For $1 \leq k \leq 6$, this was done via the approach described in Section \ref{sec:fixminindex} (hence we cannot guarantee that the computed local minimum is global); for $k\in\{0,7\}$, the problem was instead solved exactly using the singular value decompositions of the compound matrices $\begin{bmatrix}
A \\ B
\end{bmatrix}$ or $\begin{bmatrix}
A & B
\end{bmatrix}$. The results are reported in Table \ref{tab:fixminind}. Interestingly, our computations suggest that the sets $\text{Schur}_8^{3}$ and $\text{Schur}_8^{4}$ (that includes our best guess for the nearest singular pencil) are almost as close to the pencil \eqref{eq:difficultpencil}: this might explain why this example is particularly challenging.

\begin{table}[h]
\begin{center}
\begin{tabular}{c|cccc} \toprule 
    {Minimal index} &$k=0$&$k=1$&$k=2$&$k=3$\\ 
    \hline
    {Distance}  & 0.0112695  & 0.0112680 & 0.0111718&0.0111731 \\ \midrule
        {Minimal index}&$k=4$&$k=5$&$k=6$&$k=7$  \\
        \hline
    {Distance}&0.0456669&
    0.0475071&0.0477320&0.0494382\\ \bottomrule
\end{tabular}
\caption {Estimated distances of the pencil \eqref{eq:difficultpencil} from the sets $\text{Schur}_8^{k+1}=\overline{\mathcal{S}_8^k}$} \label{tab:fixminind}
\end{center}
\end{table}

In Subsection \ref{sec:multiple_starting} we saw that running the algorithm several times from different random initial points improves the solution significantly. One might wonder if we can boost the performance even further by exploiting the possibility of fixing the minimal index of the sought nearest singular pencil. To explore this question, we compare the following two different computational strategies: (S1) run the standard algorithm presented in Section~\ref{sec:complex} $n$ times starting from different starting points and choose the best solution; and (S2) run the algorithm presented in Section~\ref{sec:fixminindex} once for each minimal index and choose the best solution. We compare these two strategies on $10^3$ randomly generated $20 \times 20$ pencils, and show the results in Table \ref{tab:strategy_comparison}. 

Table \ref{tab:strategy_comparison} shows that the two strategies performed almost identically in terms of the average and median distances. Moreover, strategy (S1) performed better $50.7 \%$ of the time while strategy (S2) performed better $47.2 \%$ of the time; the two strategies computed the same minimum $2.1 \%$ of the time. For the average of $10^3$ Bernoulli distributed random variables with mean $\frac{1}{2}$, the standard deviation is approximately $1.6 \%$. Accounting for the omitted $2.1 \%$, the difference between each strategy and the mean is approximately one sigma. Consequently, the choice between these two strategies is not statistically significant.

\begin{table}[h]
\begin{center}
\begin{tabular}{cccc} \toprule
    {Strategy} & {Frequency of best output} & {Median distance} & {Average distance}  \\ \midrule
    {(S1)}  & 50.7 \%  & 1.4925 & 1.4933  \\
    {(S2)}  & 47.2 \%  & 1.4944 & 1.4964 \\ \bottomrule
\end{tabular}
\caption {Comparison of the strategies (1) and (2). The two strategies computed the same minimum $2.1 \%$ of the time. } \label{tab:strategy_comparison}
\end{center}
\end{table}

{
\subsection{Smoothed algorithm} \label{sec:smooth_numerical}
We compare the standard algorithm of Section~\ref{sec:complex} with the smoothed algorithm outlined in Section \ref{sec:smoothalgorithm}. We ran both algorithms for the same inputs, and recorded the outputs as well as the time to converge. This was done for $100$ randomly generated pencils of sizes $10$, $40$ and $80$. For the smoothed algorithm, the parameter $\alpha$ was set to $\alpha = -10^{20}$ (the code adaptively makes this value milder if the algorithm becomes numerically non-smooth; the motivation behind this choice is simply to make it the smallest possible in practice).

The average running times are shown in Table \ref{tab:smooth_comparison}. It can be seen that the difference in the running time grows larger as the size of the pencil is increased. For pencils of size of 80, the smoothed algorithm is more than $50 \%$ slower than the standard algorithm. The difference in the computed distance was not statistically significant. Moreover, the two algorithms gave the same minimizer $86.3 \%$ of the time. 
\begin{table}[h]
\begin{center}
\begin{tabular}{cccc} \toprule
    {Size} & {10} & {40} & {80}  \\ \midrule
    {Standard}  & 0.568 s  & 14.0 s & 97.9 s \\
    {Smooth}  & 0.609 s  & 18.0 s & 153.4 s \\
    {Difference}  & 7.21 \%  & 29.1 \% & 56.7 \% \\
    \bottomrule
\end{tabular}
\caption {Comparison of the running times between the standard algorithm of Section~\ref{sec:complex} and the smoothed algorithm outlined in Section \ref{sec:smoothalgorithm}. The computation was run using MATLAB R2023a on an Intel Core i5-12600K. } \label{tab:smooth_comparison}
\end{center}
\end{table}

}

\section{Conclusion}\label{sec:conclusion}

We have described a novel { paradigm} to compute the nearest singular pencil to a given one, based on Riemannian optimization. { This new approach} makes it practically feasible, for the first time, to solve the problem for pencils of moderate size, say, a few hundreds rows/columns. However, benchmarking against the top competitors \cite{bora,glm} with these sizes is computationally very demanding for \cite{bora} and not practically feasible for \cite{glm}, due to the longer running times of these methods. We have instead run an extensive statistical comparison for inputs of smaller size $n \leq 50$.  The Riemannian method does better than its competitors in terms of quality of the output when the size of the problem is not extremely small.  Furthermore, the performance is also very favourable to the new algorithm in terms of computational time: for example, on average for randomly generated inputs of size $n=50$ and using one of the authors' personal computer,  the new method converged in about $25$ seconds while \cite{bora} required more than $12$ minutes (the lack of an efficient implementation made it impractical to run \cite{glm} for $n=50$, but it already took several minutes to our implementation of \cite{glm} to solve inputs of size $n=6$, when \cite{bora} and the Riemannian method only required fractions of a second). We also demonstrated that, due to how non-convex the problem is, faster running time can be translated into better results by running the algorithm multiple times.  We have also proved that the new Riemannian approach can be adapted with minor algorithmic modifications to compute the nearest singular pencil with prescribed (right) minimal index. To this goal, we have proved several new results on the generalized Schur form of singular pencils. { One of these results also implies that the new Riemannian approach remains essentially the same for finding the nearest real singular pencil to a given real one. In fact, the general complex algorithm solves automatically such real problem when the input pencil is real and the starting point is also real.}

\section*{Acknowledgements} We sincerely thank Biswajit Das and Shreemayee Bora for providing the MATLAB codes for their algorithm in \cite{bora} and Nicola Guglielmi for detailed advice and discussions on the algorithm in \cite{glm}. {We are grateful to two anonymous reviewers whose comments led to significant improvements of our paper.} Finally, we acknowledge the computational resources provided by the Aalto Science-IT project.

\appendix

\section{``Swapping'' diagonal elements}
\label{sec:swapping}
We prove a theorem on $2 \times 2$ real or complex regular pencils that allows to interchange the position of two eigenvalues in the generalized Schur form of a \emph{regular} pencil. 

\begin{theorem} \label{th.11}
Let $A + z B$ be a $2\times 2$ upper triangular complex (resp. real) regular pencil
\[
A + z B = \begin{bmatrix}
            a_{11} + z b_{11} & a_{12} + z b_{12}  \\
            0 & a_{22} + z b_{22}
          \end{bmatrix}.
\]
Then, there exist $U,V \in \mathbb{C}^{2 \times 2}$ unitary (resp. $U,V \in \mathbb{R}^{2 \times 2}$ orthogonal) such that
\begin{equation} \label{eq.main11}
U^* (A + z B) V =  \begin{bmatrix}
            c_{11} + z d_{11} & c_{12} + z d_{12}  \\
            0 & c_{22} + z d_{22}
          \end{bmatrix} =: C + z D
\end{equation}
and
  \begin{enumerate}
    \item[(a)] $\displaystyle \frac{a_{22}}{b_{22}} = \frac{c_{11}}{d_{11}}$, that is, $a_{22} + z b_{22}$ and $c_{11} + z d_{11}$ have the same eigenvalue (which can be equal to $\infty$), and
    \item[(b)] $\displaystyle \frac{a_{11}}{b_{11}} = \frac{c_{22}}{d_{22}}$, that is, $a_{11} + z b_{11}$ and $c_{22} + z d_{22}$ have the same eigenvalue (which can be equal to $\infty$).
  \end{enumerate}
\end{theorem}

\begin{proof} The regularity of $A + zB$ implies that $(a_{11},b_{11}) \ne (0,0)$.
The matrix $S = -b_{11} A + a_{11} B$ is singular because its first column is zero. Therefore there exists $u_2 \in \mathbb{C}^{2 \times 1}$ (resp. $u_2 \in \mathbb{R}^{2 \times 1}$) with $\|u_2 \|_2 = 1$ (where $\|u\|_2 := \sqrt{u_2^* u_2}$) such that
\begin{equation} \label{eigenvector2}
0 = u_2^* S = -b_{11} (u_2^* A) + a_{11} (u_2^* B).
\end{equation}
Therefore,
\begin{align}
u_2^* A = \frac{a_{11}}{b_{11}} \, u_2^* B, & \quad \mbox{if $b_{11} \ne 0$} , \label{beta2}\\
u_2^* B = \frac{b_{11}}{a_{11}}  \, u_2^* A. & \quad \mbox{if $a_{11} \ne 0$} . \label{alpha2}
\end{align}
Let
\begin{equation}\label{U2} U = \begin{bmatrix}
           u_1 & u_2
         \end{bmatrix} \in \mathbb{C}^{2 \times 2} \quad \mbox{(resp. $\in \mathbb{R}^{2 \times 2}$)}
\end{equation}
be a unitary (resp. orthogonal) matrix whose second column is $u_2$.

Now we are going to distinguish three cases.

{\em Case 1.} $u_2^* A = u_2^* B = 0$. This case is not possible since then $u_2^* (A + z B)= 0$ for all $z$ and the pencil would be singular.

{\em Case 2.} $u_2^* A \ne 0$ ($u_2^* B$ may be zero or not). Observe that \eqref{beta2} implies that $a_{11} \ne 0$ in this case,  since $(a_{11},b_{11}) \ne (0,0)$. Let $V\in \mathbb{C}^{2 \times 2}$ (resp. $V\in \mathbb{R}^{2 \times 2}$) be a unitary (resp. orthogonal ) matrix
\begin{equation}\label{V12} V = \begin{bmatrix}
           v_1 & v_2
         \end{bmatrix} \qquad \mbox{with $\displaystyle v_2 = \frac{1}{\|(u_2^* A)^*\|_2} \, (u_2^* A)^*$.}
\end{equation}
If we take $U$ as in \eqref{U2} and $V$ as in \eqref{V12}, the second row of $U^* (A + z B) V$ is, taking into account \eqref{alpha2},
$$
(u_2^* A + z \, u_2^* B ) V = (1 + z \frac{b_{11}}{a_{11}})\, (u_2^* A) V =
(1 + z \frac{b_{11}}{a_{11}})\, \begin{bmatrix}
                               0 & *
                          \end{bmatrix},
$$
showing that the pencil in \eqref{eq.main11} is upper triangular.

{\em Case 3.} $u_2^* B \ne 0$ ($u_2^* A$ may be zero or not). Observe that \eqref{alpha2} implies that $b_{11} \ne 0$ in this case. Let $V\in \mathbb{C}^{2 \times 2}$ (resp. $V\in \mathbb{R}^{2 \times 2}$) be a unitary (resp. orthogonal) matrix
\begin{equation}\label{V22} V = \begin{bmatrix}
           v_1 & v_2
         \end{bmatrix} \qquad \mbox{with $\displaystyle v_2 = \frac{1}{\|(u_2^* B)^*\|_2} \, (u_2^* B)^*$.}
\end{equation}
If we take $U$ as in \eqref{U2} and $V$ as in \eqref{V22}, the second row of $U^* (A + z B) V$ is, taking into account \eqref{beta2},
$$
(u_2^* A + z \, u_2^* B ) V = (\frac{a_{11}}{b_{11}} + z) \, (u_2^* B) V =
(\frac{a_{11}}{b_{11}} + z) \begin{bmatrix}
                               0 & *
                          \end{bmatrix}
$$
and this proves that the pencil in \eqref{eq.main11} is upper triangular.

From \eqref{eigenvector2} and \eqref{eq.main11}, we obtain
$$
0 = u_2^* (-b_{11} A + a_{11} B) v_2 = -b_{11} c_{22} + a_{11} d_{22} \, .
$$
This proves item (b). Item (a) follows from the fact that the multiset of eigenvalues of an upper triangular regular pencil is the union of the eigenvalues of its diagonal entries.
\end{proof}

Note that one can embed the $2 \times 2$ unitary matrices of Theorem \ref{th.11} into an $n \times n$ larger unitary matrix, so that any pair of eigenvalues can be swapped. By repeated applications of this technique, it is clear that there exists a unitary equivalence that can be applied to a given particular Schur form of a pencil, so that another triangular pencil is obtained whose diagonal elements correspond to any arbitrary permutation of the eigenvalues of the pencil. We state this formally as a corollary below.

\begin{corollary}\label{cora2}
    Suppose that $A+z B$ is an $n \times n$ regular complex pencil having eigenvalues $\lambda_1,\lambda_2,\dots,\lambda_n$ (possibly with repetitions and possibly equal to $\infty$). For every permutation $\sigma$ of $\{1,\ldots ,n\}$, there exist unitary matrices $Q$ and $Z$ such that
    \begin{enumerate}
        \item $Q(A+z B)Z$ is in Schur form, and
        \item the diagonal elements $[Q(A +z B)Z]_{ii}=:t_{ii}-z s_{ii}$ satisfy
    \[ \frac{t_{ii}}{s_{ii}} = \lambda_{\sigma(i)}   \]
    for all $i=1,\dots,n$.
    \end{enumerate} 
\end{corollary}

{ 
\section{Euclidean gradient and Hessian of the smoothed objective function} \label{sec:gradient_smooth}
The $i$th entry of the gradient of the Boltzmann operator that appears in the definition of $f_\alpha (Q,Z)$ in \eqref{eq:smooth} is
% \begin{equation}
%       f(Q,Z) = \sum_{i>j} (|(QAZ)_{ij}|^2 + |(QBZ)_{ij}|^2) + \mathcal{S}_\alpha\left(x_1, \ldots, x_n\right),
% \end{equation}
% $$
% \mathcal{S}_\alpha\left(x_1, \ldots, x_n\right)=\frac{\sum_{i=1}^n x_i e^{\alpha x_i}}{\sum_{i=1}^n e^{\alpha x_i}}.
% $$
\begin{align*}
\nabla_{x_i} \mathcal{S}_\alpha\left(x_1, \ldots, x_n\right)=\frac{e^{\alpha x_i}}{\sum_{j=1}^n e^{\alpha x_j}}\left[1+\alpha\left(x_i - \mathcal{S}_\alpha\left(x_1, \ldots, x_n\right)\right)\right].
\end{align*}
 Using this expression, it can be proved that the gradient of the objective function \eqref{eq:smooth} with respect to the first argument $Q$ is
\begin{align*}
    \nabla_{Q} f_\alpha (Q,Z) &= 2 L_s (QAZ) (AZ)^* + 2 L_s (QBZ) (BZ)^* \\
    &+ 2 \mbox{Diag}(\nabla \mathcal{S}_\alpha\left(x_1, \ldots, x_n\right)) (L_d(QAZ)(AZ)^* + L_d(QBZ)(BZ)^*),
\end{align*}
where $x_i = |(QAZ)_{ii}|^2 + |(QBZ)_{ii}|^2$; the operator $L_s$ now takes the strictly lower triangular part of the matrix and sets every other element to be zero; $L_d$ takes the diagonal of the matrix and sets every other element to zero; and $\mbox{Diag}$ is the operator that turns a vector into a square diagonal matrix with the entries of the vector on the main diagonal. A similar expression holds for $\nabla_{Z} f_\alpha (Q,Z)$, and finally the Euclidean gradient of $f_\alpha (Q,Z)$ is given by 
$\nabla  f_\alpha (Q,Z) = (\nabla_Q  f_\alpha (Q,Z), \nabla_Z  f_\alpha (Q,Z)) $.

Recall from \eqref{eq:euclideanhess} that the Euclidean Hessian, considered as a linear map, is the directional derivative of the gradient and that, in our case, $d_Q$ and $d_Z$ denote the matrix directions for the matrices $Q$ and $Z$, respectively. Note that the inverse $\mbox{Diag}^{-1}$ of the $\mbox{Diag}$ operator defined above transforms a diagonal matrix into a vector given by the main diagonal of the matrix. Let us first make the following definitions:
\begin{align*}
    &\mathrm{D}\mbox{Diag}(QAZ)[d_Q,d_Z] := \mbox{Diag}^{-1}(L_d(d_Q AZ + QA d_Z)), \\
    &\mathrm{D}\mbox{Diag}(QBZ)[d_Q,d_Z] := \mbox{Diag}^{-1}(L_d(d_Q BZ + QB d_Z)),\\
    &\mathrm{D}x(Q,Z)[d_Q,d_Z] := 2 \mbox{Re}(\mbox{Diag}^{-1}(\overline{L_d(QBZ)}) \circ \mathrm{D}\mbox{Diag}(QBZ)[d_Q,d_Z] \\
    & \phantom{aaaaaaaaaaaaaaaa} + \mbox{Diag}^{-1}(\overline{L_d (QAZ)}) \circ \mathrm{D}\mbox{Diag}(QAZ)[d_Q,d_Z]),\\
    & \tilde J := (I_n + \alpha(\mbox{Diag}(x) - I_n \mathcal{S}_\alpha\left(x_1, \ldots, x_n\right))) / \sum_{j=1}^n e^{\alpha x_j}
    + (I_n - \mathbf{1} \nabla \mathcal{S}_\alpha\left(x_1, \ldots, x_n\right)^T) / \sum_{j=1}^n e^{\alpha x_j} \\
    & \quad - (\mathbf{1} + \alpha(x - \mathbf{1} \mathcal{S}_\alpha\left(x_1, \ldots, x_n\right)))\left(e^{\alpha x}\right)^T / \left(\sum_{j=1}^n e^{\alpha x_j}\right)^2, \\
    & J := \mbox{Diag}(\alpha e^{\alpha x}) \tilde J, \\
    & \mathrm{D}\nabla \mathcal{S}_\alpha[d_Q,d_Z] := J \mathrm{D}x(Q,Z)[d_Q,d_Z],
\end{align*}
where $\mathbf{1}$ denotes a column vector of ones of size $n \times 1$, $x := (x_1, \ldots, x_n)^T \in \mathbb{C}^{n \times 1}$, $e^{\alpha x} := (e^{\alpha x_1}, \ldots, e^{\alpha x_n})^T \in \mathbb{C}^{n \times 1}$ and $\circ$ is the Hadamard product. The directional derivative of $\nabla_{Q} f_\alpha (Q,Z)$ is then
\begin{align*}
    \mathrm{D} \nabla_{Q} f_\alpha (Q,Z)[d_Q,d_Z] &= 2 L_s(d_Q AZ)(AZ)^* + 2 L_s(QA d_Z)(AZ)^* + 2 L_s(QAZ)(A d_Z)^*\\
    &+ 2 L_s(d_Q BZ)(BZ)^* + 2 L_s(QB d_Z)(BZ)^* + 2 L_s (QBZ)(B d_Z)^* \\
    & + 2\mbox{Diag}(\nabla \mathcal{S}_\alpha\left(x_1, \ldots, x_n\right)) (L_d(d_QAZ)(AZ)^* + L_d(d_QBZ)(BZ)^*) \\
    & + 2\mbox{Diag}(\nabla \mathcal{S}_\alpha\left(x_1, \ldots, x_n\right)) (L_d(QA d_Z)(AZ)^* + L_d(Q B d_Z)(BZ)^*) \\ 
    &+ 2\mbox{Diag}(\nabla \mathcal{S}_\alpha\left(x_1, \ldots, x_n\right)) (L_d(QAZ)(Ad_Z)^* + L_d(QBZ)(B d_Z)^*)\\
    & + 2\mbox{Diag}(\mathrm{D}\nabla \mathcal{S}_\alpha[d_Q,d_Z]) (L_d(QAZ)(AZ)^* + L_d(QBZ)(BZ)^*).
\end{align*}
A similar result holds for the directional derivative of $\nabla_{Z} f_\alpha (Q,Z)$, which yields the Euclidean Hessian via \eqref{eq:euclideanhess}.}
\end{document}